\documentstyle[12pt,a4cpde,math,mlmak,cpde]{article}
\def\tT{\tilde T}
\def\TT0{\tilde T_0}
\def\TS{\tilde S}
\def\TA{\tilde A}
\def\TH{\tilde H}
\def\TPi{\tilde \Pi}
\def\TP{\tilde P}
\def\DST{\displaystyle}
\def\SST{\scriptstyle}

\def\htot{{H_{\rm tot}}}
\begin{document}
\thispagestyle{empty}\addtolength{\topmargin}{1.5cm}
\begin{center}{\large\bf ON INDEX FORMULAS FOR MANIFOLDS\\[0.5em]
WITH METRIC HORNS}
\end{center}

\vspace*{1.5\baselineskip}
\centerline{\large\it  Matthias Lesch and Norbert Peyerimhoff}

\vspace*{0.7\baselineskip}
\noindent
\begin{center}
\parbox[t]{6cm}{\begin{raggedright}
Institut f\"ur Mathematik\\
Humboldt--Universit\"at zu Berlin\\
Unter den Linden 6\\
D--10099 Berlin\\[0.5em]
lesch@mathematik.hu-berlin.de\\
\end{raggedright}}
\parbox[t]{5cm}{\begin{raggedright}
Ruhr-Universit\"at Bochum \\
Fakult\"at f\"ur Mathematik \\
Universit\"atsstr. 150 \\
D-44780 Bochum\\[0.5em]
peyerim@math.ruhr-uni-bochum.de
\end{raggedright}}
\end{center}

\vspace*{2\baselineskip}


\begin{abstract}
In this paper we discuss the index problem for geometric
differential operators (Spin--Dirac operator,
Gau{\ss}-Bonnet operator, Signature operator)
on manifolds with metric horns. On singular manifolds these 
operators in general do not have unique closed extensions. But there
always exist two extremal extensions $D_{\rm min}$ and $D_{\rm max}$. We
describe the quotient ${\cal D}(D_{\rm max}) / {\cal D}(D_{\rm min})$
explicitely in geometric resp. topological terms of the base manifolds of the
metric horns. We derive index formulas for the Spin--Dirac and
Gau{\ss}-Bonnet operator. For the Signature operator we present
a partial result. 
\end{abstract} 


\section{Introduction}
It is well known, that an elliptic differential operator
\[ D: C_0^\infty(M,E) \to C_0^\infty(M,F) \]
between sections of vector bundles $E$ and $F$ over a compact, closed
manifold $M$ has a unique closed extension which is a Fredholm
operator. In the more general case of an open manifold
\[ M \ = \ M_\epsilon \ \dot \cup \ U \]
\textheight21.6cm
with a compact part $M_\epsilon$ (with nonempty boundary $N$) and an open
part $U$ which we consider as a punctured neighborhood of the
``singularities'' the situation is more complicated: There may be many
closed extensions between two extremal extensions $D_{\rm min}$ and
$D_{\rm max}$. These extensions can be parametrized by closed
subspaces $V$ of the quotient ${\cal D}(D_{\rm max}) /
{\cal D}(D_{\rm min})$ ($\cal D$ denotes the domain of the operator).
It is natural to ask for a characterization of
${\cal D}(D_{\rm max}) /
{\cal D}(D_{\rm min})$ for special ``singular'' manifolds as well as for the Fredholm property and
explicit index formulas of the assigned closed extensions $D_V$ of a
given operator $D$. The complete answers to these and related questions
for the Gau{\ss}-Bonnet operator, the Signature operator and the Spin--Dirac
operator (which we refer to as the ``geometric operators'')
on manifolds with cone like singularities are given in the papers
\cite{Che0}, \cite{Che1}, \cite{Che2}, \cite{Chou}, \cite{BS}, \cite{B1}.
In this paper we consider manifolds with {\em metric horns}
in the restricted sense of \cite{Che1}. This means that $U \cong (0,\epsilon)
\times N$ with warped product metric
\[ g \ = \ dx^2 \ + \ h(x)^2 g_N, \]
and $h(x) = x^\alpha$, $\alpha > 1$. On these manifolds all the geometric
operators have a common normal form on $U$
\begin{equation} \label{nf}
D \vert U\ \cong\  \Bigl({d \over dx}\ +\ {h'(x) \over h(x)} S\Bigr)\ 
\oplus \ \Bigl( {d \over dx}\ +\ {1 \over h(x)} \tilde S\ +\
{h'(x) \over h(x)} \tilde A(x)\Bigr).
\end{equation}
This will be discussed in detail in Section 2.
Section 3 investigates the nature of ${\cal D}(D_{\rm max}) /
{\cal D}(D_{\rm min})$ and the Fredholm
property. In Section 4 we elaborate on the index problem.
The key tool is a 
homotopy transforming the operator of interest
into a so called regular singular operator. Regular singular operators
turn up in the conic situation, for which index formulas exist
(cf. the papers mentioned above). With the help of a result in \cite{CL} we
will show that the index remains constant during the homotopy.
By this procedure we bypass the need of an explicit functional calculus
on metric horns, which seems to be very difficult
and not yet developed to the knowledge of the authors. (Cf. also the remark in
\cite[pp.138]{Che1}: {\it ``At the level of an explicit functional calculus
the more difficult case of metric horns must be distinguished from the
special case of metric cones...''}) The last two sections 
are devoted to the geometric operators as concrete examples.

\addtolength{\topmargin}{-1.5cm}
We point out as a remarkable fact that the situation of metric horns
turns out to be much simpler than the case of metric cones.
For example, by looking at the space ${\cal D}(D_{\rm max}) /
{\cal D}(D_{\rm min})$ describing the variety of
closed extensions, we prove the following characterizations
for the geometric operators on manifolds $M$ with metric horns:\pagebreak[4]
\[ \widearray{\cal D}(D_{\rm max}) /
{\cal D}(D_{\rm min}) \ \cong \ \left\{ \begin{array}{cl} \{ 0 \} & \mbox{for the Signature
operator},\\ \{ 0 \} & \mbox{for the Gau{\ss}-Bonnet operator ($m$ even)},\\
{\cal H}^{n/2}(N) & \mbox{for the Gau{\ss}-Bonnet operator ($m$ odd)},\\
{\rm ker}\ D_N & \mbox{for the Spin--Dirac operator ($m$ even),}
\end{array} \right.\]
where $m={\rm dim}\ M$, $n={\rm dim}\ N$ and $D_N$ denotes the Spin--Dirac
operator on the base manifold $N$ of the horns. Let us emphasize, that
the nature of ${\cal D}(D_{\rm max}) /
{\cal D}(D_{\rm min})$ depends only on the topology of $N$ for the first
two operators. It is an interesting question, whether there exists a canonical 1:1
correspondence between ideal boundary conditions and closed extensions
of the Gau{\ss}-Bonnet operator in this particular situation 
(cf. \cite[Theorem 3.8 a)]{BL2}.

For the Spin--Dirac
operator ${\rm ker}\ D_N$ is not a topological invariant (cf.
\cite{H}). However, we have at least independence of ${\rm dim}\ 
{\cal D}(D_{\rm max}) /
{\cal D}(D_{\rm min})$ under
conformal changes of the metric on $N$ for this operator.
This yields a remarkable contrast to the cone like case, where always small
eigenvalues (in the interval $(-1/2,1/2)$) are involved. They are obviously
not stable under even such an easy deformation as multiplying the metric
of $N$ by a constant factor. 

In Section 5 we derive an index formula for the Spin--Dirac operator:
$$ {\rm ind}\ D^+_{\rm max/min} \ = \ \int_M \hat A \ - \ {\eta(0)
\over 2} \ \pm \ {b \over 2}$$
with $b:=\dim\ker D_N$.

The index
formula for the
Gau{\ss}-Bonnet operator on manifolds with horns developed in Section 6
becomes simply
\[ {\rm ind}\ D_{\rm GB} \ = \ \int_M e \ + \ {1 \over 2} \Bigl(
\sum_{k=0}^{{m\over 2}-1} (-1)^k {\rm dim}\ H^k(N) \ - \
\sum_{k={m\over 2}}^n (-1)^k {\rm dim}\ H^k(N)  \Bigr). \]

We emphasize that this formula is not based on our analysis but
is derived by reinterpretation of results of J. Cheeger.
In contrast to the cone like case there occurs no additional
spectral data of $N$ and we conclude in particular $\int_M e \in {\Bbb{Z}}$.

At this point we would like to make some historical remarks:
A detailed analytic study (Hodge theory,
functional calculus, index formulas, ...) of inductively defined singular
spaces based on manifolds with metric cones was first carried out by J.
Cheeger in the papers \cite{Che0,Che1,Che2}. In the same spirit
\cite{Chou} considered
the Spin--Dirac operator, stressing the importance of small
eigenvalues of the operator on the base manifold of the cones.
The paper \cite{BS} derived an index formula for differential
operators with a special normal form.
The framework of \cite{BS} covers the geometric operators on manifolds
with asymptotic conic singularities and motivated the present
paper to a great extend. In the meantime
the regular singular case was generalized 
in \cite{B1} and \cite{B2} and became applicable
to singular algebraic curves (cf. \cite{BPS})
and to manifolds with finitely many ends
and some further geometric conditions. Another generalization of the regular
singular case to operators of Fuchsian type is part of \cite{L1}.  
Homotopy arguments based on the result of \cite[Theorem 4.1.]{CL} and
similar to those in this article  
have been carried out in \cite{SSi}, \cite{S} and \cite{B2}.

\medskip

Part of this article (concerning the Spin--Dirac
operator on metric horns) is one of the main results of the thesis of
the second named author.

\bigskip
\noindent
{\large\bf List of Notations}

\medskip
\noindent
\begin{tabular}{llp{9.8cm}}
$\cd(D)$ & domain of the (unbounded) operator $D$\\
$D_\min, D_\max$ & minimal/maximal closed extension of a differential
 operator\\
 $D^t$ & formal adjoint of the differential operator $D$\\
 $\cl(H)$ & algebra of bounded operators on the Hilbert space $H$\\
 $X$ & operator of multiplication by $x$\\
 $\Rightarrow$& uniform convergence\\
 $\|\cdot\|_{\rm HS}$ & Hilbert--Schmidt norm\\
 $\sigma(S)$& spectrum of the self--adjoint operator $S$\\
 $P^*$ & Hilbert space adjoint of the (possibly unbounded) operator P
\end{tabular}



\section{A generalization of regular singular operators}
We consider the following situation:
\begin{quote}
Let $(M,g_M)$ be a (singular) Riemannian manifold of dimension $m$, $E,F$
two Riemannian resp. hermitian vector bundles over $M$ and $D:
C^\infty_0(E)
\to C^\infty_0(F)$ a first order elliptic differential operator. Let
$U$ be a neighborhood of the singularities, such that $M \backslash U$
is compact with smooth compact boundary $N$. Furthermore, let $G$ be a
Riemannian resp. hermitian vector bundle over $N$, such that we obtain
the following identifications by a suitable unitary separation of
variables (cf. \cite{BS}):
$$\widearray\begin{array}{rcll}
\DST U &\cong& (0,\epsilon) \times N \qquad& {\rm with}\ 0 < \epsilon \le
1,\\
\DST D &\cong& {d \over dx} + {\alpha \over x} S & {\rm on\ the\ open\ set}\
U.
  \end{array}$$
$S: C^\infty(G) \to C^\infty(G)$ is supposed to be a first order symmetric
elliptic differential operator on $N$. By a slight abuse of notation, there
exists an orthonormal basis of $L^2(G)$, $\{ e_s \}$, of eigenfunctions of $S$.
\end{quote}
At this point the meaning of the constant $\alpha \ge 1$ is not apparent
because we could write $S$ instead of $\alpha S$. However, the parameter
$\alpha$ will play a crucial role in the case of a metric horn. To avoid
notational confusion, it is introduced here.

This is the easiest
example of a so called {\em regular singular operator}. More
general operators of regular singular type are treated in
\cite{BS}, \cite{B1}, \cite{B2} or \cite{L1}. The example above serves as a
starting point for
a general framework, which applies to the ``geometric operators''
(i.e. Spin--Dirac operator, Gau{\ss}-Bonnet operator and Signature operator)
on manifolds with metric horns. However, in order to treat the case of horns
we will have to leave the regular singular situation.

Before doing so, let us introduce some 
results of \cite{BS} and \cite{B1} which we shall need:


\begin{enumerate}
\item Let $H_{1 \over 2\alpha} :=
\bigoplus_{\vert s \vert < 1/(2\alpha)} \ker(S - s)$ and consider the Hilbert spaces
${\cal D} (D_{\rm max})$ and ${\cal D}(D^t_{\rm max})$
equipped with graph-norms. Then the linear functionals
\begin{eqnarray*}
&\Phi:& \left\{ \renewcommand{\arraystretch}{1.5} \begin{array}{ccc}
{\cal D}(D_{\rm max}) & \to & H_{1\over 2\alpha} \\
f & \mapsto & \displaystyle{\bigoplus_{\vert s \vert < 1/(2\alpha)} 
\lim_{x \to 0} x^{\alpha s}
\ \Pr_{\ker(S-s)} f(x)}, \end{array} \renewcommand{\arraystretch}{1}
\right. \\
&\Phi':& \left\{ \renewcommand{\arraystretch}{1.5} \begin{array}{ccc}
{\cal D}(D^t_{\rm max}) & \to & H_{1\over 2\alpha} \\
f & \mapsto & \displaystyle{\bigoplus_{\vert s \vert < 1/(2\alpha)} 
\lim_{x \to 0} x^{-\alpha s}
\ \Pr_{\ker(S-s)} f(x)}, \end{array} \renewcommand{\arraystretch}{1}
\right.
\end{eqnarray*}
are well defined and continuous. Here $\Pr_{W}: L^2(G) \to
W$ denotes the orthogonal projection onto $W \subset H_{1 \over 2\alpha}$.
\item For all $f \in {\cal D}(D_{\rm max})$ and $g \in {\cal D}(D^t_{\rm max})$
\[ ( D_{\rm max} f , g ) \ = \ ( f , D^t_{\rm
max} g )\ -\ ( \Phi f , \Phi' g )_{H_{1 \over 2\alpha}}. \]
\item All closed
extensions of $D$ are Fredholm. These extensions are in 1:1 cor\-res\-pondence
to subspaces of $H_{1 \over 2\alpha}$. The corresponding extension 
to a subspace 
$W \subset H_{1 \over 2\alpha}$ is given by
\[ {\cal D}(D_W) \ := \ \{ f \in {\cal D}(D_{\rm max}) \ \vert \ \Phi f
\in W \}. \]
Any orthogonal decomposition $H_{1 \over 2\alpha} = V \perp W$ 
implies: $D^t_V = (D_W)^*$.
\end{enumerate}

{\it We now describe the generalizations which will be applicable to metric
horns.} We consider again the situation described at the beginning of
this section, but now with a more general
\begin{equation}
D: C_0^\infty(E) \to C_0^\infty(F).
\label{MLrev1}
\end{equation} 
We assume that after a suitable
choice of isometries of Hilbert spaces
\begin{equation}
   \begin{array}{rcl}
   \DST L^2(E|U)&\cong &\DST L^2((0,\epsilon),L^2(G)),\\[0.5em]
   \DST L^2(F|U)&\cong &\DST L^2((0,\epsilon),L^2(G))
   \end{array}
  \label{MLrev2}
\end{equation}
and a suitable splitting $L^2(G)=H_{\rm tot} = 
H\oplus\tilde H$ ($H, \tilde H$ Hilbert
spaces), $D|U$ decomposes into $T \oplus \tilde T$ with
\begin{equation}\begin{array}{rcl}
  \DST T &=& \DST {d \over dx}\ +\ {h'(x) \over h(x)} S,\\[1em]
\DST \tilde T &=& \DST {d \over dx}\ +\ {1 \over h(x)} \tilde S\ +\
{h'(x) \over h(x)} \tilde A(x),
  \end{array}
  \label{MLrev3}
\end{equation}
and $h(x) = x^\alpha$ for $\alpha > 1$.
The operator $T: C^\infty_0((0,\epsilon),H)$ is
the same as introduced at the beginning since ${h' \over h} = {\alpha
\over x}$. So
$S$ is supposed to be a self-adjoint operator on ${\cal D}(S)$
with an appropriate family of eigenfunctions $e_s$. In the case of
geometric differential operators on metric horns it turns
out that $H$ is of finite dimension which makes the regular
singular component of the operator $D|U$ particularly easy to handle.

Let us turn over to the properties of $\tilde T$. For notational
convenience we put a tilde $\tilde {\ }$ over all notions connected with
this operator. $\tilde S$ is a symmetric operator defined on a dense
subspace ${\cal D}(\tilde S)$ of $\tilde H$. As before, let $\{ \tilde e_s \}$
be an orthonormal basis of $\tilde H$ consisting of eigenfunctions of
$\tilde S$. Additionally, we assume ${\rm ker}\ \tilde S = \{ 0 \}$.
$\tilde A(\cdot)$ is a smooth family of bounded operators on $\tilde H$:
$\tilde A \in C^\infty( (0,\epsilon), {\cal L}(\tilde H) )$. Moreover,
there exists a constant $C > 0$ with
\[ \Vert \tilde A(x) \Vert \ \le \ C \quad {\rm for\ all}\ x \in
(0,\epsilon). \]
In the sequel we consider the term ${h' \over h} \tilde A$ as a
perturbation and denote by
\begin{equation}
\tilde T_0 \ = \ {d \over dx} + {1 \over h(x)} \tilde S \label{MLrev3b}
\end{equation}
the ``unperturbed part of $\tilde T$''. This point of view turns out to
make sense as long as $\alpha > 1$.

The geometric operators on manifolds with metric horns
are examples of this general framework. Let us discuss this in more
detail:

\begin{description}
\item[Spin--Dirac operator] Let $M$ be an open spin manifold of even dimension
with metric $g_M$, where $U$ and $N$ are as in the beginning of this
section. Assume the
existence of a fixed metric $g_N$ with
\begin{equation}\label{metrics}
g_M \vert U \ \cong \ dx^2\ + \ h(x)^2 g_N.
\end{equation}
Denote by $S(M) \to M$ an irreducible spin bundle and $D: C_0^\infty(S(M))
\to C_0^\infty(S(M))$ the corresponding Spin--Dirac operator. As it
is well known, $D$ splits into
\[ D \ = \ \left( \begin{array}{cc} 0 & D^- \\ D^+ & 0 \end{array}
\right) \]
with respect to the eigenbundles $S^+(M), S^-(M)$ of the multiplication
with the complex volume element. By a suitable
unitary separation of variables $D^+$ can be transformed into
(see \cite[p. 652]{S})
\begin{equation}\label{sdo}
D^+ \ \cong \ {d \over dx} \ + \ {1 \over h(x)} D_N \qquad {\rm on}\
U,
\end{equation}
where $D_N: C^\infty(S(N)) \to C^\infty(S(N))$ denotes the Spin--Dirac
operator on $N$ with the induced spin structure. (We use the
orientation of \cite{APS} which is opposite to that in \cite{S}.)
Consequently, we choose
$H := {\rm ker}\ D_N$, $\tilde H$ its orthogonal complement in
$L^2(S(N))$ and $S = 0$, $\tilde S := D_N \vert \tilde H$, $\tilde
A = 0$.
\item[Gau{\ss}-Bonnet operator] This operator is given by
\[ D_{\rm GB} \ := \ d_M + d_M^t: \ \Omega_0^{\rm even}(M) \ \to \
\Omega_0^{\rm odd}(M). \]
Let $M$ be an open Riemannian manifold and $U, N, g_M$ and $g_N$ etc. 
as above. The calculations in \cite[pp. 696]{BS} yield:
\[ D_{\rm GB} \vert U \ \cong \ T_{\rm GB} : \ C^\infty((0,\epsilon),
\bigoplus_{j=0}^n \Omega^j(N)) \ \to \ C^\infty((0,\epsilon),
\bigoplus_{j=0}^n \Omega^j(N)) \]
with $n = {\rm dim}\ N$, $T_{\rm GB} = {d \over dx} + {h'(x) \over h(x)}
S_1 + {1 \over h(x)} S_2$ and
\begin{eqnarray*}
S_1 &=& {\rm diag}(c_j)_{0 \le j \le n},\quad c_j \ := \ (-1)^j
(j-{n \over 2}),\\
S_2 &=& d_N \ + \ d_N^t \ = \ \left( \begin{array}{cccc}
0 & d_N^t & & \\ d_N & \ddots & \ddots & \\ & \ddots & \ddots
& d_N^t \\ & & d_N & 0 \end{array} \right).
\end{eqnarray*}
Obviously, $S_2$ is a square root of the Laplacian
on forms. (For metric collars, i.e. $h = 1$,
the situation simplifies to $D_{\em GB} \vert U \cong {d \over dx}
+ S_2$.)

In the case of metric horns, i.e. $h(x) = x^\alpha$, we choose
$H$ to be the space ${\cal H}(N)$ of harmonic forms and $\tilde H$ its
orthogonal complement in $L^2(\Lambda^* T^* N)$. $H$ and $\tilde H$ are
obviously invariant subspaces of $S_1$ and $S_2$ and $S_2 \vert H
= 0$. Moreover, we choose
\[ S\ := \ S_1 \vert H, \quad \tilde S \ := \ S_2 \vert \tilde H,
\quad \tilde A \ := \ S_1 \vert \tilde H. \]

\item[Signature operator] We assume $M$ to be an oriented, open
Riemannian manifold of dimension $m=4k$, $U$ and $N$ as above and that
\myref{metrics} holds. The involution
\[ \tau: \ \left\{ \begin{array}{c} \Omega(M) \ \to \ \Omega(M), \\
\tau \alpha \ := \ i^{2k+j(j+1)} *_M \alpha \quad {\rm for}\ \alpha \in
\Omega^j(M) \end{array} \right. \]
anti-commutes with the operator $d_M + d_M^t: \Omega_0(M) \to \Omega_0(M)$.
Hence $d_M + d_M^t$ interchanges the $\pm 1$-eigenspaces
$\Omega_0^\pm(M)$ of $\tau$. The restriction
\[ D_{\rm S}:\ d_M + d_M^t:\ \Omega_0^+(M) \ \to \ \Omega_0^-(M) \]
is called the Signature operator. Unitary separation of variables
yields according to \cite[pp.707]{BS}:
\[ D_{\rm S} \vert U \ \cong \ T_S \ = \ {d \over dx}\ + \
{h'(x) \over h(x)} S_1\ + \ {1 \over h(x)} S_2 \]
with $S_1 = {\rm diag} (b_j)_{0 \le j \le n}$, $n= {\rm dim}\ N$,
$b_j = {n \over 2} - j$ and
\[ S_2 \vert \Omega^j(N) \ = \ (-1)^{k+1+[ {j+1 \over 2} ]} \bigl(
(-1)^j *_N d_N - d_N *_N \bigr).  \]
$S_2$ is symmetric and $S_2^2 = \Delta_N$. Consequently, $H := {\cal
H}(N)$ and $\tilde H := {\cal H}(N)^\perp$ are invariant subspaces of
$S_1$ and $S_2$. Choose
\[ S\ := \ S_1 \vert H, \quad \tilde S \ := \ S_2\vert \tilde H, \quad
\tilde A \ := \ S_1 \vert \tilde H. \]
\end{description}

It was pointed out by the referee that the metric horn case can be
'reduced' to the conic case by a logarithmic change of variables.
Namely, given the metric
\[ g \ = \ dx^2 \ + \ x^{2\alpha}g_N, \quad \alpha > 1.\]
Putting $x(t) = \bigl( (1-\alpha)\log t \bigr)^{1 \over 1-\alpha}$ we
find
\begin{eqnarray*}
g &=& {x(t)^{2\alpha} \over t^2} \bigl( dt^2 \ + \ t^2 g_N \bigr)\\
  &=:& \rho(t)^2 \bigl( dt^2 \ + \ t^2 g_N \bigr),
\end{eqnarray*}
which is a conic metric, however with a singular conformal factor.
Applying the standard unitary separation of variables (cf. \cite[p.443]{BL2})
to the Gau{\ss}-Bonnet operator, one finds
\[ D_{GB} \ = \ {1 \over \rho}\ ( {d \over dt} \ + \ {1 \over t} (S_1 + S_2)
\ + \ {\rho' \over \rho} S_3 ), \]
where $S_1, S_2$ are as in the above description of the Gau{\ss}-Bonnet 
operator and $S_3 = {\rm diag}(\tilde c_j)_{0\le j \le n}$, $\tilde c_j :=
(-1)^j (j-{n \over 2}) - {1 \over 2}$. However, ${1 / \rho}$, ${\rho'
/ \rho^2}$ are singular at $0$ and the usual regular singular analysis
cannot be applied directly. 

We turn back to our general D described in \myref{MLrev1}--\myref{MLrev3}.
In Section 3 we will characterize the $L^2$--closed 
extensions of $D$ and show that these are Fredholm. Section 4 deals with a
particular family $\{ D_\beta \}_{\beta \in [ \beta_1,\beta_2 ] }$ of operators
and we will show that the index of corresponding $L^2$-closed extensions
remains constant under variation of $\beta$. 
In both sections we need to construct parametrices and to prove certain properties about them. The construction is done by
patching together an interior pseudo--differential parametrix
and a suitable boundary parametrix. The boundary parametrix will
be obtained by considering $D|U$ as an infinite sum of one--dimensional
differential operators (cf. \myref{MLrev5} below). 
For the construction of the boundary parametrix
we have to introduce some integral operators:

Similar to \cite[4.5/4.6]{LSE}, 
to a given function $F \in C^\infty(0,1)$ and $s \in {\Bbb R}$ 
we introduce (whenever these integrals exist)
\begin{equation} \label{parametrix}
\widearray\begin{array}{rcl}
\DST P_{0,s}^F f(x) &:=& \DST \int_0^x \exp\Bigl( - \int_y^x s F(t) dt \Bigr)
f(y) dy, \\[1em]
\DST P_{1,s}^F f(x) &:=& \DST \int_1^x \exp\bigl( - \int_y^x s F(t) dt \Bigr)
f(y) dy,
\end{array}
\end{equation}
for all $f \in L^2(0,1)$. 

In the particular case $F(x)=1/x$ these operators become
\[ \begin{array}{rcll}
   \DST P^F_{0,s} f (x) &:=& \DST\int_0^x \bigl( {y \over x} \bigr)^s f(y) dy
     &{\rm for} \ s > -{1 \over 2}, \\[1em]
   \DST P^F_{1,s} f (x) &:=& \DST \int_1^x \bigl( {y \over x} \bigr)^s f(y) dy
     &{\rm for \ all} \ s.
     \end{array}
\]
These operators were already used in \cite[(2.3), (2.4)]{BS}.

One easily checks
\begin{equation} \label{TPeqI}\label{MLrev4}
   \widearray\begin{array}{l}
   \DST \bigl( {d \over dx} + F s \bigr)\ P_{j,s}^F f = f,\quad\mbox{for}\quad
    f\in L^2(0,1),\\
   \DST  P_{j,s}^F \bigl( {d \over dx} + F s \bigr) f = f,\quad\mbox{for}\quad
    f\in C_0^\infty(0,1),\\
    \DST \bigl(P_{0,s}^F\bigr)^* = - P_{1,-s}^F.
    \end{array}
\end{equation}

We introduce the unperturbed operators
\begin{equation}
   \TT0:= \frac {d}{dx}+\frac{1}{h(x)}\TS,\quad D_0:=T\oplus \tT_0,
   \label{ML2}
\end{equation}
where $D_0$ lives on $U$.

Furthermore, we denote by $\Pi, \TPi$ the orthogonal projections from
$\htot:=H\oplus \TH$ onto $H$ resp. $\TH$.

Via the spectral decomposition of $S$, $\tilde S$
the operators $T, \tilde T_0$ can be viewed as the infinite direct sum
\begin{equation}
    \moplus_{s\in\spec S} \left(\frac{d}{dx} + \frac{h'}{h}s\right)
    \oplus \moplus_{s\in\spec\tilde S}\left(\frac{d}{dx}+\frac{1}{h} s\right).
\label{MLrev5}
\end{equation}
Thus it seems natural to construct parametrices for $T, \tT$ as follows: 
Put
\begin{eqnarray}
    P:=P(T)&:=
    & \Bigl( \DST\moplus_{\begin{array}{c}\SST s< s_1\\
      \SST s\in\spec(S)
      \end{array}} P^{h'/h}_{1,s} \Bigr)\ \moplus \ \Bigl( 
\moplus_{\begin{array}{c}\SST s> s_1\\
       \SST s\in\spec(S)\end{array}} P^{h'/h}_{0,s} \Bigr),
    \label{MLrev6}
\end{eqnarray}   
where $s_1\in [-1/(2\alpha),1/(2\alpha)]\setminus \spec (S)$. The choice of
$s_1$ corresponds to a particular choice of a $L^2$--closed
extension of $D$.
Furthermore, the parametrix for $\tilde T$ is:
\begin{equation}
  \TP:=\TP(\tT):=\Bigl( \moplus_{\begin{array}{c}\SST s< 0\\
      \SST s\in\spec( \tilde S)
      \end{array}} P_{1,s}^{1/h} \Bigr)
             \ \oplus \ \Bigl( \moplus_{\begin{array}{c}\SST s>0\\
      \SST s\in\spec( \tilde S)
      \end{array}} P_{0,s}^{1/h} \Bigr).
  \label{ML7} \label{NPrev2}
\end{equation}
The boundary parametrix is then given by $P_{\rm bd} = P \oplus \tilde P$
and all operators $P^{h'/h}_{j,s}, P_{j,s}^{1/h}$ in this construction 
are Hilbert-Schmidt. 
We will use the following notation for particular choices of $s_1$ in
\myref{MLrev6}:

\begin{equation}\begin{array}{rcl}
\DST   P^\max:=P^\max(T)&:=&\Bigl(\DST\moplus_{\begin{array}{c}\SST s< 
1/(2\alpha)\\
      \SST s\in\spec(S)
      \end{array}} P^{h'/h}_{1,s}\Bigr)\ \oplus \ \Bigl(
\moplus_{\begin{array}{c}\SST s\ge 1/(2\alpha)\\
      \SST s\in\spec(S)
      \end{array}} P^{h'/h}_{0,s} \Bigr),\\[3em]
\DST   P^\min:=P^\min(T)&:=&\Bigl(\DST\moplus_{\begin{array}{c}\SST s\le 
-1/(2\alpha)\\
      \SST s\in\spec(S)
      \end{array}} P^{h'/h}_{1,s}\Bigr)\ \oplus \ \Bigl( 
\moplus_{\begin{array}{c}\SST s> -1/(2\alpha)\\
      \SST s\in\spec(S)
      \end{array}} P^{h'/h}_{0,s} \Bigr),\\[3em]
\DST   P_\delta:=P_\delta(T)&:=&\Bigl( 
  \DST\moplus_{\begin{array}{c}\SST s<0\\
      \SST s\in\spec(S)
      \end{array}} P^{h'/h}_{1,s}\Bigr)\ \oplus \ \Bigl(
\moplus_{\begin{array}{c}\SST s\ge 0\\
      \SST s\in\spec(S)
      \end{array}} P^{h'/h}_{0,s}\Bigr).
   \end{array}\label{ML6}
\end{equation}
In view of \myref{MLrev4} we have
\begin{equation}
P^{\max/\min}(T)^*=-P^{\min/\max}(-T^t), \quad \TP(\tT)^*=-\TP(-\tT^t).
  \label{ML8}
\end{equation}


\section{The $L^2$--closed extensions of $D$}

The goal of this section is to classify the $L^2$-closed
extensions of $D: C^\infty_0(E) \to C^\infty_0(F)$ and to prove their
Fredholm property. In doing so, we follow step by step the scheme of
\cite{BS}. It turns out, that the calculations for the second component
$\tilde T$ can be carried out in a very similar way as in the
regular singular case. Most of the calculations are even more simple
than in the regular singular case.

In this section we relax the axioms \myref{MLrev3}  
for $D$ slightly. Namely, we assume that $D \vert U$ decomposes into
$T \oplus \tilde T$ with
\begin{eqnarray*}
T &=& {d \over dx} \ + \ {\alpha \over x} S,\\
\tilde T &=& {d \over dx} \ + \ {1 \over h(x)} \tilde S \ + \ 
{1 \over x} \tilde A(x),
\end{eqnarray*}
where $h$ is a positive 
increasing function with $h(x)\le x^\ga$ for some $\ga>1$.
Moreover, we assume that
$(0,1)\ni x\mapsto (I+|\TS|)^{-1}\TA(x)$ and $(0,1)\ni x\mapsto \TA(x)(I+|\TS|)^{-1}$
are smooth maps into $\cl(\TH)$ and that
\begin{equation}
    \sup_{0<x<1}\|(I+|\TS|)^{-1}\TA(x)\|+\|\TA(x)(I+|\TS|)^{-1}\|\le C.
    \label{lastlabel}
\end{equation}    

We will use the abbreviation
\begin{equation}
\mu(x):=\int_x^1 \frac{1}{h(y)} dy.
 \label{ML3}
\end{equation} 
Note that 
\begin{equation}
    \lim_{x\to 0+} \mu(x)=+\infty.
  \label{MLrev9}
\end{equation}

\begin{lemma} \label{MLS1} 
{\rm (cf. \cite[Lemma 4.1]{LSE})} For $f\in L^2(0,1)$ we have
the following estimates:
\renewcommand{\labelenumi}{{\rm (\roman{enumi})}}
\begin{enumerate}
\item $\displaystyle |P_{0,s}^{1/h} f(x)|\le \frac{1}{\sqrt{2 s}}\,
          \sqrt{h(x)}\, \|f\|_{L^2}, \quad s>0,$
\item $\displaystyle |P_{1,s}^{1/h} f(x)|\le \frac{c(h,s_0)}{\sqrt{|s|}}
       \, x^{\ga/2}\, \|f\|_{L^2},\quad s\le s_0<0.$
\end{enumerate}
\end{lemma}

\proof We use the Cauchy--Schwarz to obtain the estimate
\begin{equation}
|P_{i,s}^{1/h} f(x)| \le e^{s \mu(x)} \Big|\int_j^xe^{-2 s \mu(y)} dy
   \Big|^{1/2} \|f\|_{L^2},\quad j=0,1.
    \label{MLrev8}
\end{equation}
Because of the monotonicity of $h$ this implies for $j=0$
\begin{eqnarray*}
    |P_{0,s}^{1/h} f(x)| &\le& \sqrt{h(x)} e^{s \mu(x)} \Big(\int_0^x
     \frac{1}{h(x)}e^{-2 s \mu(y)} dy\Big)^{1/2} \|f\|_{L^2}\\
     &=& \frac{\sqrt{h(x)}}{\sqrt{2s}}\|f\|_{L^2},
\end{eqnarray*}
and we have proved the first inequality.

In view of \myref{MLrev8} to prove (ii) it suffices to prove
$$e^{2 s \mu(x)} \int_x^1 e^{-2s \mu(y)} dy \le \frac{c(h,s_0)}{|s|} x^\ga.$$
To do this we split the integral:
\begin{eqnarray*}
   e^{2s \mu(x)} \int_x^{2x} e^{-2s \mu(y)} dy &=& e^{2s \mu(x)}\int_x^{2x}
     h(y) \frac{1}{h(y)} e^{-2s \mu(y)} dy\\
     &\le& \frac{h(2x)}{2 |s|} \big(1- e^{2|s| (\mu(2x)-\mu(x))}\big)\\
     &\le& \frac{h(2x)}{2|s|}\le \frac{2^{\ga-1}}{|s|} x^\ga.
\end{eqnarray*}

Similarly we find the estimate
\begin{eqnarray*}
     e^{2s \mu(x)} \int_{2x}^1 e^{-2s \mu(y)} dy
     &\le& \frac{h(1)}{2 |s|} \big(e^{2 |s| (\mu(2x)-\mu(x))}- e^{2|s| \mu(x)}\big)\\
     &\le& \frac{c(h)}{2|s|}e^{-2|s| \int_x^{2x} \frac{dy}{h(y)}}\\
     &\le& \frac{c(h,s_0)}{|s|} x^\ga.
\end{eqnarray*}
Here we have used the inequality
$$e^{-|s| \int_x^{2x} \frac{dy}{h(y)}}\le c(h,s_0,\gamma) x^\gamma,\quad
   |s|\ge s_0$$
for any $\gamma\ge 0$. This inequality is an easy consequence of
$h(x)\le x^\ga$.\endproof

\begin{lemma}
   \label{MLrev-Lemma15}
   Let $\gb>-\ga$ and $0<\epsilon\le 1$. Then in $L^2(0,\epsilon)$ 
we have
$$ \begin{array}{rclcll}
    \DST \|X^\gb P_{0,s}^{1/h}\|&+ &\| P_{1,-s}^{1/h}X^\gb\|&\le&\DST \frac 1s 
               \eps^{\ga+\gb},& s>0,\\
    \DST  \|X^\gb P_{1,s}^{1/h}\|&+& \|P_{0,-s}^{1/h}X^\gb\|&\le&\DST 
          \frac{c(h,\gb,s_0)}{|s|} \eps^{\ga+\gb},& s\le s_0<0.
   \end{array}$$
\end{lemma}
\proof Since $(X^\gb P_{0,s}^{1/h})^*= P_{1,-s}^{1/h}X^\gb,
(X^\gb P_{1,s}^{1/h})^*= P_{0,-s}^{1/h}X^\gb,$ it suffices to
estimate $\|X^\gb P_{0,s}^{1/h}\|$ for $s>0$ and
$\|P_{1,s}^{1/h}X^\gb\|$ for $s\ge |s_0|>0$. Since the operators
are integral operators with non--negative kernels, we may apply
Schur's test \cite[p. 22]{HS}.

The kernel of $X^\gb P_{0,s}^{1/h}$ is given by
$$k_1(x,y)=x^\gb e^{s(\mu(x)-\mu(y))},\quad x>y.$$
We have
\begin{eqnarray*}
   \int_0^\eps k_1(x,y) dy &\le& x^\gb e^{s \mu(x)} \int_0^x y^\ga
     \frac{1}{h(y)} e^{-s\mu(y)} dy\\
    &\le & x^{\ga+\gb} \frac 1s \le \frac{\eps^{\ga+\gb}}{s},
\end{eqnarray*}
and
\begin{eqnarray*}
    \int_0^\eps k_1(x,y) dx &\le& e^{-s \mu(y)} \int_y^\eps x^{\ga+\gb}
     \frac{1}{h(x)} e^{s\mu(x)} dx\\
    &\le & \frac{\eps^{\ga+\gb}}{s} \big(1-e^{s(\mu(\eps)-\mu(y))}\big)\\
   &\le& \frac{\eps^{\ga+\gb}}{s},
\end{eqnarray*}
hence Schur's test implies $\|X^\gb P_{0,s}^{1/h}\|\le 
\frac{\eps^{\ga+\gb}}{s}.$

Next let $s\ge |s_0|>0$ and consider the kernel of $P_{0,s}^{1/h}X^{\gb}$,
$$k_2(x,y)=y^\gb e^{s(\mu(x)-\mu(y))},\quad x>y.$$
Then similarly as before
$$\int_0^\eps k_2(x,y) dy \le \frac{\eps^{\ga+\gb}}{s},$$
and similar as in the proof of Lemma \ref{MLS1} we estimate
\begin{eqnarray*}
    \int_0^\eps k_2(x,y) dx &\le& 
      y^\gb e^{-s \mu(y)}\int_y^\eps x^\ga
     \frac{1}{h(x)} e^{s\mu(x)} dx\\
    &\le & y^{\gb} \frac {e^{-s\mu(y)}}{s}
    \left( 2^\ga y^\ga(e^{s\mu(y)}-e^{s\mu(2y)})+
      \eps^\ga (e^{s\mu(2y)}-e^{s\mu(\eps)})\right)\\
    &\le &\frac{c(h,\gb,s_0)}{s} \eps^{\ga+\gb},
\end{eqnarray*}
and invoking again Schur's test we reach the conclusion.
\endproof

Now we introduce the operators
\begin{equation}
P_{\rm bd}^{\max/\min}:= P^{\max/\min}\oplus \TP.
  \label{ML9}
\end{equation}  

For the next lemma recall the notation $\htot = H \oplus \tilde H$. $\Pi$ and
$\tilde \Pi$ are the corresponding orthogonal projections.

\begin{lemma} \label{MLS2} {\rm 1.} For $f\in L^2((0,1),\htot)$ we have the estimates
\begin{equation}\widearray\begin{array}{rcl}
 \DST  \| \Pi (I+|S|)^{1/2} P_{\rm bd}^\min f(x)\|&\le& \DST c\, \sqrt{x}\, \|f\|,\\
 \DST  \| \tilde\Pi (I+|\TS|)^{1/2} P_{\rm bd}^\min f(x)\|&\le& \DST c\, x^{\ga/2}\, \|f\|.
 \end{array}\label{ML10}
\end{equation}

{\rm 2.} For $0< \eps\le 1$ we have
$$\| X^{-1} \TA P_{\rm bd}^{\min/\max}\|_{L^2((0,\eps),\htot)}+
  \|P_{\rm bd}^{\min/\max} X^{-1} \TA \|_{L^2((0,\eps),\htot)}\le c \,\eps^{(\ga-1)}.$$

{\rm 3. } $P_{\rm bd}^\max$ maps $L^2((0,1), \htot)$ into $\cd(D_{0,\max})$.
\end{lemma}

\proof 1. The first inequality follows from \cite[Lemma 2.1]{BS}, the second is an immediate
consequence of Lemma \ref{MLS1}.

2. This is a consequence of Lemma \ref{MLrev-Lemma15} and
\myref{lastlabel}.

3. This is an immediate consequence of \myref{TPeqI} and Lemma
\ref{MLS1}.
\endproof

\begin{lemma} \label{MLS3} {\rm 1. (cf. \cite[Lemma 3.6]{B1})} Let $\varphi\in \cinfz{[0,1)}$
with $\varphi=1$ near $0$. Then we have for $f\in \cd(D_\max)$
$$P_{\rm bd}^\max D_\max \varphi f = \varphi f + (P_{\rm bd}^\max X^{-1} \TA)\varphi f.$$

{\rm 2. (cf. \cite[Lemma 3.4]{B1})} There exists an $\eps>0$, such that
for $\varphi, \omega\in\cinfz{[0,\eps)}$, $\varphi, \omega=1$ near
$0$ and $\omega \varphi=\varphi$, we have for $f\in \cd(D_\max)$
$$\varphi f= \omega P_{\rm bd}^\max V D_\max \varphi f$$
with some bounded operator $V$.

\end{lemma}

\proof 1. Since the projections $\Pi, \TPi$ commute with $P_{\rm bd}^\max, D_\max$
it follows from \cite[Sec. 3]{B1} that
$$\Pi P_{\rm bd}^\max D_\max \varphi f=\varphi \Pi f.$$
Now we put
$$g:= \varphi f + P_{\rm bd}^\max X^{-1} \TA \varphi f$$
and $g_s(x)=(g(x), \tilde e_s)$ where $\tilde e_s$ is an eigenvector of $\TS$,
$\TS \tilde e_s=s \tilde e_s$. We have $g_s\in L^2(0,1)$ and, moreover,
$$ g_s'(x) +\frac{s}{h(x)} g_s(x) = ((D_\max \varphi f)(x), \tilde e_s).$$
Thus $g_s$ and $((P_{\rm bd}^\max D_\max \varphi f)(x), \tilde e_s)$ satisfy the same
first order differential equation. If $s<0$, then
$$g_s(1)= \Big(P_{1,s}^{1/h}(X^{-1} \TA \varphi f (\cdot),
  \tilde e_s)\Big)(1)=0$$
and
$$((P_{\rm bd}^\max D_\max \varphi f)(1), \tilde e_s)=
   P_{1, s}^{1/h}(D_\max \varphi f, \tilde e_s)(1)=0$$
hence we have equality in this case.

If $s>0$, then the $g_s$ and
$((P_{\rm bd}^\max D_\max \varphi f)(x), \tilde e_s)$ differ by $c\, 
e^{s \mu}$, which is square integrable iff $c=0$. Since $g_s$ and
$((P_{\rm bd}^\max D_\max \varphi f)(x), \tilde e_s)$ are square integrable they
must be equal.

2. In view of Lemma \ref{MLS2}, 2. we choose $\eps>0$ such that
$$\| X^{-1} \TA P_{\rm bd}^{\min/\max}\|_{L^2((0,\eps),\htot)}+
  \|P_{\rm bd}^{\min/\max} X^{-1} \TA \|_{L^2((0,\eps),\htot)}< \frac 12.$$
Furthermore, let
$\omega, \psi, \varphi\in\cinfz{[0,\eps}$, $\omega=\psi=\varphi=1$
near $0$ and $\omega\psi=\psi, \psi\varphi=\varphi$.

By 1.
\begin{eqnarray*}
   (\omega P_{\rm bd}^\max D_\max \psi)\varphi f&=& \varphi f + (\omega P_{\rm bd}^\max X^{-1}
    \TA \psi)\varphi f\\
    &=:&(I+R) \varphi f
\end{eqnarray*}
and $\|R\|<1/2$.

Introducing $R_1:= \psi X^{-1} \TA P_{\rm bd}^\max \omega$, we also have
$\|R_1\|<1/2$. By induction, one easily finds
$$R^n (\omega P_{\rm bd}^\max D_\max \psi)\varphi f= \omega P_{\rm bd}^\max R_1^n D_\max \varphi f.$$

Thus we conclude
\begin{eqnarray*}
    \varphi f&=& (I+R)^{-1}(\omega P_{\rm bd}^\max D_\max \psi)\varphi f\\
        &=&\sum_{n=0}^\infty (-1)^n R^n (\omega P_{\rm bd}^\max D_\max \psi)\varphi f\\
        &=&\omega P_{\rm bd}^\max D_\max \varphi f+
           \omega P_{\rm bd}^\max \sum_{n=1}^\infty (-1)^n R_1^n D_\max\varphi f\\
       &=:& \omega P_{\rm bd}^\max V D_\max \varphi f
\end{eqnarray*}
and we are done.

\endproof

\begin{cor} \label{MLS4} Let $\varphi$ be as before.

{\rm 1.} We have
 $$ \varphi \cd(D_{\max/\min})=\varphi \cd (D_{0,\max/\min}).$$

{\rm 2.} For $f\in \cd(D_\max)$ we have $\|\TPi f(x) \|= O(x^{\ga/2}),\; x\to 0$.

{\rm 3.} The maps $\Phi$ and $\Phi'$ are well--defined 
on $\cd(D_\max)$ and $\cd(D_\max^t)$ and for
$f\in \cd(D_\max), g\in \cd(D_\max^t)$ we have
$$(D_\max \varphi f, \varphi g)- (\varphi f, D_\max^t\varphi g)= 
(\Phi f, \Phi' g)_{H_{1 \over 2\alpha}}.$$
\end{cor}

\proof 1. If $f\in \cd(D_\max)$, then Lemma \ref{MLS3}, 2. shows that
$\varphi f\in \Im \varphi P_{\rm bd}^\max\subset \cd(D_{0,\max})$
by Lemma \ref{MLS2}, 2. On the other hand, if $f\in \cd(D_{0,\max})$,
then we have from Lemma \ref{MLS3}, 1. $P_{\rm bd}^\max D_{0,\max} \varphi f=\varphi f$ and
Lemma \ref{MLS2} implies $\varphi f \in \cd(D_\max)$.

On $\cinfz{(0,\eps),\TH}$ we have, for $\eps$ small enough
\begin{eqnarray*}
   \| X^{-1} \TA f\|&=& \| X^{-1} \TA P_{\rm bd}^\min D_{0,\min} f\|\\
     &\le & \| X^{-1} \TA P_{\rm bd}^\min \| \|D_{0,\min} f\|\\
     &<& \frac 12 \| D_{0,\min} f\|,
\end{eqnarray*}
hence $X^{-1} \TA$ is $D_0$--bounded with bound $<1/2$ which easily
implies $\varphi\cd(D_\min)=\varphi\cd(D_{0,\min})$.

2. This is an immediate consequence of Lemma \ref{MLS3} and Lemma \ref{MLS2}.

3. In view of 1. and 2. it suffices to prove this for $D_0$.

We use the decomposition $\varphi f = f_0 + \tilde f_0$ and $\varphi
g = g_0 + \tilde g_0$. As a consequence of Corollary \ref{MLS4}, 2. 
\begin{eqnarray}
( \tilde T \tilde f_0 , \tilde g_0 ) \ - \
( \tilde f_0 , \tilde T^t \tilde g_0 ) &=&
\int_0^\epsilon \partial_x \ ( \tilde f_0(x), \tilde g_0(x) )_{\tilde H}
\ dx \nonumber \\
&=& - \lim_{x \to 0} \ ( \tilde f_0(x), \tilde g_0(x) )_{\tilde H}
\nonumber \\
&=& 0.
\label{g150}
\end{eqnarray}

This implies together with result 2 of Section 2 and the identity $D^t_{\rm max} =
{D_{\rm min}}^*$
$$
( \tilde T \tilde f_0 , \tilde g_0 ) \ + \
( T f_0 , g_0 )
= ( \varphi f , D^t_{\rm max} \varphi g ) \
- \ ( \Phi f ,
\Phi' g )_{H_{1 \over 2\alpha}}. \eqno\mbox{\qed}$$

\begin{theorem} \label{class} The closed extensions of $D$ are in
$1:1$ correspon\-den\-ce
to subspaces of $H_{1\over 2\alpha}$. For a given subspace 
$W \subset H_{1 \over 2\alpha}$ the
corresponding extension is given by
\[ {\cal D}(D_W)\ :=\ \{ f \in {\cal D}(D_{\rm max})\ \vert\ \Phi f \in
W \}. \]
Any orthogonal decomposition $H_{1\over 2\alpha} = V \perp W$ 
implies $D^t_V = (D_W)^*$.
\end{theorem}

\proof Obviously, every closed extension $\bar D$ corresponds to
a closed vector space $\tilde W = {\cal D}(\bar D)$ between ${\cal
D}(D_{\rm min})$ and ${\cal D}(D_{\rm max})$ with respect to the graph
norm. So we only have to show, that
\[ \Phi: {\cal D}(D_{\rm max}) \ \to \ H_{1\over 2\alpha} \]
is a continuous epimorphism with kernel ${\cal D}(D_{\rm min})$. To
prove surjectivity we consider $f(x) := \varphi(x) x^{-\alpha s} e_s$ for any
$s \in \sigma(S)$ with $\vert s\vert < {1 \over 2\alpha}$. $f$ extends
trivially to a section $f \in L^2(E)$ with $D_{\rm max}f = \varphi' 
x^{-\alpha s} e_s \in L^2(F)$ and $\Phi f = e_s$.

Corollary \ref{MLS4}, 3. implies immediately ${\rm ker}\ \Phi 
= {\cal D}(D_{\rm min})$ and $D^t_V = {D_W}^*$ for any orthogonal 
decomposition $H_{1 \over 2\alpha} = V \perp W$.\endproof

\remark An immediate consequence of the preceding proof is
\[ {\cal D}(D_{\rm max})/{\cal D}({D_{\rm min}}) \ \cong \ 
H_{1 \over 2\alpha}. \]

\begin{theorem} \label{fred}
All $L^2$-closed extensions of $D$ are Fredholm and
\[ {\rm ind}\ D_W \ = \ {\rm ind}\ D_{\rm min} \ + \ {\rm dim}\ W \]
for all subspaces $W \subset H_{1 \over 2\alpha}$.
\end{theorem}
\proof  First, we construct a right parametrix. We choose
$\omega, \psi$ as before with $\omega \psi=\psi$, such that
for $R:=\omega X^{-1} \TA P_{\rm bd}^\min \psi$ we have
$\|R\|<\frac 12$.

Since $D$ is an elliptic operator, there exists an interior parametrix,
$P_i$, such that
$$ D P_i= I-\psi+ K$$
and $K$ is a compact pseudodifferential operator with compact support.
We put
$$Q:= P_i +\omega P_{\rm bd}^\min \psi$$
and find
\begin{eqnarray*}
  DQ&=& I-\psi +K + \omega'P_{\rm bd}^\min\psi +\omega DP_{\rm bd}^\min \psi\\
    &=& I+R+K+\omega'P_{\rm bd}^\min\psi.
\end{eqnarray*}
Now, $\omega'P_{\rm bd}^\min\psi$ is a continuous operator with
$\Im \omega'P_{\rm bd}^\min\psi\subset H_{\rm comp}^1(M,E)$, the Sobolev space
with compact support in $M$. By the closed graph theorem,
$$\omega'P_{\rm bd}^\min\psi:L^2(M,F)\to H_{\rm comp}^1(M,E)$$
is continuous and
hence $\omega'P_{\rm bd}^\min\psi$ is a compact operator $L^2(M,F)\to L^2(M,E)$.
Finally, $Q(I+R)^{-1}$ is a right parametrix for $D_{\max/\min}$.

In the same way we find a right parametrix for $D_{\max/\min}^t$. The
adjoint of this parametrix will then be a left parametrix for $D_{\max/\min}$
and we obtain the Fredholm property of $D_{\max/\min}$.

Since the inclusion map $\iota: {\cal D}(D_{\rm
min}) \to {\cal D}(D_W)$ is Fredholm with index $- {\rm dim}\ W$, we
obtain by the logarithmic law
$$ {\rm ind}\ D_{\rm min} \ = \ {\rm ind}\ D_W \ + \ {\rm ind}\ \iota \ =
\ {\rm ind}\ D_W \ - \ {\rm dim}\ W. \eqno\mbox{\qed}$$


\section{Index considerations}
In this section we assume $D: C_0^\infty(E) \to C_0^\infty(F)$ to be of
the generalized form described in Section 2 with the difference,
that $h(x) = x^\alpha$ $(\alpha > 1)$ only on a subinterval
$(0,\epsilon_0)$ of $(0,\epsilon)$. 

We restrict
ourselves to the particular extension $D_\delta := D_W$ (see Theorem 
\ref{class}) corresponding to the subspace
\[ W \ = \ \bigoplus_{-{1\over 2\alpha} < s < 0} {\rm ker}\ (S-s) \]
and denote this extension by $D_\delta$. Since we consider in this section
only this particular extension we will often omit the index $\delta$ and simply
write $D$. Furthermore we will use the notations
\[ U_{\epsilon} \ := \ (0,\epsilon) \times N \quad {\rm and} \quad 
M_{\epsilon} \ := \ M \ - \ \bigl( (0,\epsilon] \times N \bigr). \]
The following proposition
states that we can remove the perturbation ${h' \over h} \tilde A$ in the
operator of consideration without changing the index.
\begin{prop} \label{removeperturbation}
Let $\alpha > 1$ and
$\psi_0 \in C^\infty(0,\epsilon)$ with $\psi_0(\epsilon) =0$ and
$\psi_0 \vert (0,\epsilon_0) = 1$. Let
$D_0$ be an ``unperturbed version of $D$'', i.e. $D_0$ coincides with $D$
on $M_\epsilon$ and
\[ D_0 \vert U \ \cong \ \Bigl( {d \over dx} + {h' \over h} S \Bigr) \
\oplus \ \Bigl( {d \over dx} + {1 \over h} \tilde S + (1 - \psi_0) {h'
\over h} \tilde A \Bigr). \]
Then ${\rm ind}\ D_\delta = {\rm ind}\ D_{0,\delta}$.
\end{prop}
\proof From Section 3 we conclude ${\cal
D}(D_\delta)
= {\cal D}(D_{0,\delta})$ and $\tilde A_{\psi_0} := {\psi_0} {h' \over
h}\tilde A
\in {\cal L}({\cal D}(D_\delta),L^2(F))$. Here we think of $\tilde
A_{\psi_0}$ as
an operator on the whole manifold $M$ by using
the unitary separation of variables and trivial extension to $M$.
Consequently
\[ [0,1] \ni \mu \ \to \ D_{0,\delta} \ + \ \mu \tilde A_{\psi_0} \ \in
\ {\cal L}({\cal D}(D_\delta),L^2(F)) \]
is a continuous curve in the space of Fredholm operators connecting
$D_{0,\delta}$ with $D_\delta$.\endproof

Next we consider a family $\{ D^\beta \}_{\beta \in [\beta_1,\beta_2]}$
of unperturbed operators coinciding on $M_\epsilon$ and
\begin{equation} \label{operators}
D^\beta \vert U \ 
\cong \ \Bigl( {d \over dx} + {h_\beta' \over h_\beta} S
\Bigr)\ \ \oplus \ \Bigl( {d \over dx} + {1 \over h_\beta} \tilde S
\Bigr)
\end{equation}
with $h_\beta \vert (0,\epsilon_0) = x^\beta$ and $1 \le \beta_1 <
\beta_2$. Later, in our examples we will consider a particular
geometric operator on a fixed manifold $M$ with continuously changing
metrics ${g_M}^\beta$ on $U$.
The resulting
identifications \myref{operators} will be based on different separations
of variables for each $\beta \in [\beta_1,\beta_2]$. After these
identifications we are independent of metric considerations on $U$ and
can work entirely in the fixed Hilbert spaces
\begin{equation} \label{Hspaces} \label{NPrev1}
\begin{array}{rcl}
\DST {\cal H}_1 &:=& L^2\bigl( (0,\epsilon), \htot \bigr) \
\oplus \ L^2(E\vert M_\epsilon), \\[1em]
\DST {\cal H}_2 &:=& L^2\bigl( (0,\epsilon), \htot \bigr) \
\oplus \ L^2(F\vert M_\epsilon)
\end{array}
\end{equation}
with $H_{\rm tot} = H \oplus \tilde H$. 
Note, that we do not exclude the case $\beta_1 =1$. In this case,
however, both components $T^{\beta_1}$ and $\tilde T_0^{\beta_1}$ 
(see \myref{MLrev3} and \myref{MLrev3b} for the definition) are
regular singular and we choose $D^{\beta_1}_\delta$ to be the
extension corresponding to the subspace
\[ W_1 \ := \ \bigoplus_{-{1 \over 2} < s < 0} {\rm ker} \ \bigl(
S - s \bigr) \ \oplus \ \bigoplus_{-{1 \over 2} < s < 0} {\rm ker}
\ \bigl( \tilde S - s \bigr). \]
This convention coincides with the operator $D_\delta$ introduced in
\cite{BS}.

Our goal is to prove ${\rm ind}\ D^{\beta_1} = {\rm ind}\ D^{\beta_2}$
under weak additional conditions. Let $F_\beta := {h_\beta ' \over
h_\beta}$ and $\tilde F_\beta := {1 \over h_\beta}$ for notational
convenience. Then
\[ D^\beta \vert U \ \cong \ T^\beta \oplus \tilde T_0^\beta \ = \ 
\Bigl( {d \over dx} + F_\beta S \Bigr)
\ \oplus \ \Bigl( {d \over dx} + \tilde F_\beta \tilde S \Bigr). \]
We assume the existence of a constant $C > 0$ with
\begin{equation} \label{cond}
F_\beta(x), \tilde F_\beta(x) \ \ge \ C \quad {\rm for\ all}\ \beta \in
[ \beta_1, \beta_2 ] \ {\rm and}\ x \in (0,\epsilon)
\end{equation}
and uniform convergence $F_\beta \Rightarrow F_\gamma$ and $\tilde
F_\beta \Rightarrow \tilde F_\gamma$ on compact intervals $[x_1,x_2]
\subset (0,\epsilon)$ for $\beta \to \gamma \in [\beta_1,\beta_2]$.
These are the only conditions we need for the proof,
that ${\rm ind}\ D^\beta$ is independent of $\beta$.

The corresponding boundary parametrices $P_{{\rm bd},\beta}$ are given by
\begin{equation}
  P_{{\rm bd},\beta} \ := \ P_\delta(T^\beta) \ \moplus \ \tilde P(\tilde
T^\beta), \label{gg105}
\end{equation}   
using the notation \myref{NPrev2} and \myref{ML6}.

Let $0 < \epsilon_0 < \epsilon_1 < \epsilon$ and 
$\varphi, \psi \in C^\infty(0,\epsilon)$ 
with $\varphi + \psi = 1$. Moreover, $\varphi(x) = 0$ for $x$ near
$0$ and $\psi(x) = 0$ for $x$ near $\epsilon$. Let
$\varphi_2, \psi_2 \in C^\infty(0,\epsilon)$ with the same properties
near $x=0$ and $x=\epsilon$ and $\varphi_2|[0,\epsilon_1]=0,$ 
$\varphi \varphi_2 = \varphi$, $\psi \psi_2
= \psi$. We extend these functions in an obvious manner to the manifold $M$.

Analogously to \cite[Lemma 2.3. and (2.12)]{BS}
we conclude
\[ {\rm im}\ \psi_2 P_{{\rm bd},\beta} \psi \ \subset \ {\cal D}(D^\beta_\delta)
\]
and
\[ D^\beta \bigl( \psi_2 P_{{\rm bd},\beta} (\psi f) \bigr) \ = \ \psi_2'
P_{{\rm bd},\beta} (\psi f) \ + \ \psi f. \]
Let $P_i: {\cal H}_2 \to H^1_{\rm loc}(E)$ be an
interior parametrix, i.e.
$$\widearray\begin{array}{rcll}
    \varphi_2 P_i ( \varphi D^\beta f ) &=& \varphi f \ + \ L_i f \quad
                      &{\rm for\ all}\ f \in {\cal D}(D^\beta), \\
     D^\beta ( \varphi_2 P_i (\varphi g) ) &=& \varphi g \ + \ R_i g &
       {\rm for\ all}\ g \in {\cal H}_2.
     \end{array}
     $$
$P_i$ does not dependent on $\beta$ and
$L_i, R_i$ as well as their adjoints are assumed to be infinitely
smoothing and compact. Obviously $L_i f\vert U_{\epsilon_1} = 0$
for any $f \in {\cal H}_1$. The same holds true for $R_i, {L_i}^*,
{R_i}^*$. The previous remarks together with the explicit
description of the adjoints in Theorem \ref{class} and
\myref{operators} imply easily (see \myref{NPrev1} for the notation ${\cal H}_1,
{\cal H}_2$)
\begin{lemma} \label{DQ}
Let $Q_\beta := \psi_2 P_{{\rm bd},\beta} \psi + \varphi_2 P_i \varphi$. 
$Q_\beta$ resp. ${Q_\beta}^*$ are maps into
${\cal D}\bigl( D^\beta \bigr)$ resp. ${\cal D}\bigl( (D^\beta)^*
\bigr)$ and
\begin{eqnarray*}
D^\beta \bigl( Q_\beta f \bigr) &=& f \ + \ \psi_2' P_{{\rm bd},\beta} (\psi f)
\ + \ R_i f \qquad {\rm for}\ f \in {\cal H}_2, \\
(D^\beta)^* \bigl( {Q_\beta}^* g \bigr) &=& g + \psi' {P_{{\rm bd},\beta}}^*
(\psi_2 g) \ + \ {L_i}^* g \qquad {\rm for}\ g \in {\cal H}_1.
\end{eqnarray*}
\end{lemma}

To prove compactness of $Q^\beta$ we will use the following special case
of Lemma \ref{MLrev-Lemma15}: 
\begin{lemma} \label{normp} Assume $F(x) \ge {1 \over C_0}$ on
$(0,\epsilon)$ and $F(x) \ge {c \over x}$ near $x=0$ for some $c>0$.
Then
\[ \Vert P_{j,s}^F \Vert \le {C_0 \over s} \qquad {\rm for \ all}
\ \ (-1)^j\cdot s > 0.\]
\end{lemma}
\begin{prop} \label{compact}
$Q^\beta$ and $\bigl( Q^\beta \bigr)^*$ are compact operators.
\end{prop}
\proof We prove compactness of $\varphi_2 P_i \varphi$ and $\psi_2
P_{{\rm bd},\beta} \psi$ separately.
The ellipticity of $D^\beta$ implies ${\cal D}\bigl( D^\beta \bigr)
\subset H_{\rm loc}^1(E)$, so compactness of $\varphi_2 P_i \varphi$
follows from the Lemma of Rellich.

>From \myref{cond} we conclude with Lemma \ref{normp}
\begin{equation}
\label{gg108}
\Vert P_{j,s}^{F_\beta} \Vert,\ \ \Vert P_{j,s}^{\tilde F_\beta}
\Vert \ \to \
0 \qquad {\rm as}\ \ s \to (-1)^j \infty.
\end{equation}
Let temporarily ${\rm Pr}_t$ denote the orthogonal projection of 
$\htot = H \oplus \tilde H$ onto the space spanned by all eigenfunctions $e_s$ 
and $\tilde e_s$ with $\vert s \vert < t$. Since the operators $P_{j,s}^{F_\beta}, P_{j,s}^{\tilde F_\beta}$ are Hilbert-Schmidt, 
$\psi_2 \bigl( P_{{\rm
bd},\beta} \circ {\rm Pr}_t \bigr) \psi$ is compact for all $t > 0$.
\myref{gg108} implies that this operators converge to $\psi_2
P_{{\rm bd},\beta} \psi$ as $t \to \infty$ in the operator norm. This proves
compactness of $\psi_2 P_{{\rm bd},\beta} \psi$.\endproof

Next we introduce the graphs
\[ {\cal G}(D^\beta) \ := \ \{ (f,D^\beta f) \ \vert \ f \in
{\cal D}(D^\beta) \} \ \subset \ {\cal H}_1 \oplus {\cal H}_2\]
and the following metric $d$ between two closed
subspaces $V$ and $W$ of ${\cal H}_1 \oplus {\cal H}_2$:
\[ d(V,W) \ = \ \sup_{x \in V,\ \Vert x \Vert = 1} d(x,W) \ + \
\sup_{y \in W,\ \Vert y \Vert = 1} d(y,V) \]
with $d(x,W) := \inf_{y \in W} \Vert x - y \Vert_{{\cal H}_1 \oplus
{\cal H}_2}$. This metric induces a topology on the closed subspaces.
Theorem 4.1. of \cite{CL} implies, that ${\rm ind}\ D^\beta$ remains
constant whenever the mapping
\begin{equation}
\label{gg125}
[\beta_1,\beta_2] \ni \beta \ \mapsto \ {\cal G}(D^\beta)
\end{equation}
is continuous.

The next lemma is exactly the note at the end of Section 2 in \cite{SSi}
which was given without a proof. For the sake of completeness, we give
a short outline of the proof.
\begin{lemma} \label{equiv} We define
\[ E_\beta \ := \ \left( \begin{array}{cc} I & -\bigl( D^\beta \bigr)^* \\
D^\beta & I \end{array} \right): \ {\cal D} \bigl( D^\beta \bigr) \oplus
{\cal D} \bigl( (D^\beta)^* \bigr) \ \to \ {\cal H}_1 \oplus {\cal H}_2. \]
$E_\beta^{-1}$ is a bounded operator, and continuity of {\myref{gg125}}
is equivalent to the continuity of
\[ [\beta_1,\beta_2] \ni \beta \ \mapsto \ E_\beta^{-1} \in {\cal L} \bigl(
{\cal H}_1 \oplus {\cal H}_2 \bigr). \]
\end{lemma}

\proof Obviously, $E_\beta$ is closed and $E_\beta \ge I$. The same
holds for ${E_\beta}^*$. From this we conclude
bijectivity of $E_\beta$. Hence $E_\beta^{-1}$ is bounded by the closed
graph theorem.

According to \cite{CL}, ${\cal G}(D^\beta) \to
{\cal G}(D^\gamma)$ is equivalent to ${\cal G}\bigl((D^\beta)^*\bigr) \to
{\cal G}\bigl((D^\gamma)^*\bigr)$.
The equivalence to ${\cal G}(E_\beta) \to
{\cal G}(E_\gamma)$ follows from explicit formulas for the
metric, i.e. the fact that
\[ \Vert (f,\tilde f) - (g,\tilde g) \Vert \ + \
\Vert \bigl( f - {D^\beta}^* \tilde f, D^\beta f - \tilde f \bigr) -
\bigl( g - {D^\gamma}^* \tilde g, D^\gamma g - \tilde g \bigr) \Vert \]
can be estimated from below and above by multiples of
\[ \Vert f - g \Vert \ + \ \Vert D^\beta f - D^\gamma g \Vert \ + \
\Vert \tilde f - \tilde g \Vert \ + \ \Vert (D^\gamma)^* \tilde g -
(D^\beta)^* \tilde f \Vert. \]
According to \cite{CL}, ${\cal G}(E_\beta) \to {\cal G}(E_\gamma)$
is equivalent to ${\cal G}(E_\beta^{-1}) \to
{\cal G}(E_\gamma^{-1})$, which finally is equivalent to the
convergence in norm $E_\beta^{-1} \to E_\gamma^{-1}$ by the addendum
of \cite{CL}. \endproof

The uniform convergence of $F_\beta$ and $\tilde F_\beta$
on compact intervals is needed in the following
lemma which is in some sense the heart of the proof of
${\rm ind}\ D^{\beta_1} \ = \ {\rm ind}\ D^{\beta_2}$.
\begin{lemma} \label{pcont} The mapping
\[ [\beta_1,\beta_2] \ni \beta \mapsto P_{{\rm bd},\beta} 
\in {\cal L}\bigl( L^2( (0,\epsilon), \htot ) \bigr) \]
is continuous.
\end{lemma}
\proof It is enough to prove that, for each $\delta_0
> 0$, there exists a $\mu > 0$ with
\begin{equation}
\label{ggg10}
\Vert P_{0,s}^{F_\beta} - P_{0,s}^{F_\gamma} \Vert_{\rm HS}^2 < \delta_0
\ \ {\rm and} \ \ \Vert P_{0,s}^{\tilde F_\beta} - P_{0,s}^{\tilde
F_\gamma} \Vert_{\rm HS}^2  < \delta_0
\end{equation}
for all $\vert \beta - \gamma \vert < \mu$. The corresponding statement
for $P_{1,s}^{F_\beta}, P_{1,s}^{\tilde F_\beta}$ follows by taking adjoints. The proof proceeds in
two steps. First we prove convergence of the Hilbert-Schmidt norms as $\beta
\to \gamma$ for each $s$-value individually.
A {\em contraction property} allows us in the
second step to conclude uniform convergence for all $s$-values by considering
only finitely many $s$-values.

\smallskip

\noindent
{\bf (i)} The uniform convergence on compact intervals
implies for $\beta \to \gamma$
\begin{equation}\label{lala}
k_\beta(x,y) := {\rm exp}\bigl( -\int_y^x s F_\beta(t)\
dt \bigr) \ \longrightarrow \
k_\gamma(x,y) := {\rm exp}\bigl( -\int_y^x s F_\gamma(t)
\ dt \bigr)
\end{equation}
pointwise for all $0 < y \le x < \epsilon$. From $F_\beta \ge 0$ we
conclude 
\[ \vert k_\beta(x,y) - k_\gamma(x,y) \vert \ \le \ 1, \]
and thus by Lebesgue's dominated convergence theorem
\[ \Vert P_{0,s}^{F_\beta} - P_{0,s}^{F_\gamma}
\Vert_{\rm HS}^2 \ \to 0. \]
Similarly, we obtain $\Vert P_{0,s}^{\tilde F_\beta} - 
P_{0,s}^{\tilde F_\gamma} \Vert_{\rm HS} \to 0$.

\smallskip

\noindent
{\bf (ii)} Let $\gamma \in [\beta_1,\beta_2]$ and $\delta_0 > 0$ be
given. Choose a function $F \in C^\infty\bigl( (0,\epsilon) \bigr)$
with $F(x) \ge {1 \over x}$ near $x=0$, $F$ being bounded away from $0$
and
\[ F \le F_\beta \qquad {\rm for\ all}\ \beta \in [\beta_1,\beta_2]. \]
Let $k(x,y) := \exp \bigl( -\int_y^x F(t) dt \bigr)$ and $k_\beta(x,y)
:= \exp \bigl( -\int_y^x F_\beta(t) dt \bigr)$. (Differently to
\myref{lala} this time there is no `s' in the definition of $k_\beta$!)

\begin{figure}[here]
\unitlength1in
\begin{picture}(5,3)
\put(1,0.75){\vector(1,0){3}}
\put(4.05,0.7){$x$}
\put(1.5,0.25){\vector(0,1){2.5}}
\put(1.45,2.8){$y$}
\put(3,0.7){\line(0,1){0.1}}
\put(2.95,0.6){$\epsilon$}
\put(1.45,2.25){\line(1,0){0.1}}
\put(1.3,2.2){$\epsilon$}
\thicklines
\linethickness{0.5mm}
\put(1.5,0.75){\line(1,0){1.5}}
\put(1.5,0.75){\line(1,1){1.5}}
\thinlines
\put(1.75,0.75){\vector(0,1){0.25}}
\put(2,0.75){\vector(0,1){0.5}}
\put(2.25,0.75){\vector(0,1){0.75}}
\put(2.5,0.75){\vector(0,1){1.0}}
\put(2.75,0.75){\vector(0,1){1.25}}
\put(3,0.75){\vector(0,1){1.5}}
\put(1.85,0.4){$\scriptstyle \tilde k \equiv 0$}
\put(2,0.5){\vector(1,1){0.2}}
\put(1.85,1.75){$\scriptstyle \tilde k \equiv 1$}
\put(2.1,1.8){\vector(1,-1){0.2}}
\put(3.2,1.2){\shortstack[l]{$\scriptstyle \tilde k$ \scriptsize
strictly increasing \\ \scriptsize along the arrows}}
\put(4,0.2){\bf figure 1}
\end{picture}
\end{figure}


Figure 1 illustrates the behaviour of $k$ over the set $\Delta :=
\{ (x,y)\ \vert \ 0 \le y \le$ $x \le \epsilon \}$. Since
$k$
is continuous on $\Delta - \{(0,0)\}$, there exists a $w > 0$ with
\[ {\rm vol} \Bigl( A(w) := \{ (x,y) \ \vert \ k(x,y) > 1 - w \} \Bigr)
\ < \ {\delta_0 \over 2}. \]
Furthermore, there exists a $s_w > 0$, such that for all $s \ge s_w$ the
mappings
\[ \Lambda_s: \left\{ \begin{array}{ccc} [0,1-w] & \to & [0,1-w], \\
x & \mapsto & x^s \end{array} \right. \]
are contractions. From this we conclude
\begin{eqnarray*}
\lefteqn{\Vert P_{0,s}^{F_\beta} -  P_{0,s}^{F_\gamma} \Vert_{\rm HS}^2 \ 
 = \ \int_\Delta \underbrace{\vert \Lambda_s(k_\beta(x,y)) -
\Lambda_s(k_\gamma(x,y)) \vert^2}_{\le 1} d(x,y)} \\
& & \hspace{1cm} \le \ \int_{\Delta - A(w)} \vert k_\beta(x,y) -
k_\gamma(x,y) \vert^2 \ d(x,y) \ + \ {\rm vol} \Bigl( A(w) \Bigr)
\\
& & \hspace{1cm} < \ \Vert P_{0,1}^{F_\beta} - P_{0,1}^{F_\gamma} 
\Vert_{\rm HS}^2 \ + \ {\delta_0 \over 2}
\end{eqnarray*}
for all $s \ge s_w$.
By (i) we can obviously find a $\mu$ with $\Vert P_{0,1}^{F_\beta}
 -  P_{0,1}^{F_\gamma} \Vert_{\rm HS}^2 < {\delta_0 \over 2}$
for $\vert \beta - \gamma \vert < \mu$. This implies
$\Vert P_{0,s}^{F_\beta} - P_{0,s}^{F_\gamma} \Vert_{\rm HS}^2 <
\delta_0$
for all $\vert \beta - \gamma \vert < \mu$ and $s \ge
s_w$. Applying (i) again to the finitely many $s$-values below $s_w$
finishes the proof. The same holds true for the functions $\tilde
F_\beta$. \endproof

The remaining arguments follow exactly the scheme of
\cite[pp.285]{B2}. To keep the notation of loc. cit., we denote by $F^\beta$
the following operator
\[ F^\beta := \left( \begin{array}{cc} 0 &
Q_\beta \\ - {Q_\beta}^* & 0 \end{array} \right). \]
(To avoid notational misunderstandings let us remark, that the functions
$F_\beta$ are no longer used in this section.) $Z^\beta$ denotes the
orthogonal projection onto ${\rm ker}\ F^\beta$. Then
\[ E^\beta F^\beta \ = \ I + F^\beta + G^\beta \quad {\rm with}\
G^\beta := \left( \begin{array}{cc} {L_i}^* + \psi' P_{{\rm bd},\beta}^*
\psi_2 & 0 \\ 0 & R_i + \psi_2' P_{{\rm bd},\beta} 
\psi \end{array} \right). \]
The compactness of $G^\beta$ follows easily, if we replace $\psi_2
P_{{\rm bd},\beta} \psi$ by $\psi_2' P_{{\rm bd},\beta} \psi$ 
resp. $\psi_2 P_{{\rm bd},\beta} \psi'$ in the proof of 
Proposition \ref{compact}.
\begin{lemma} \label{propz} $Z^\beta$ has the following properties:
\begin{enumerate}\renewcommand{\labelenumi}{{\rm (\roman{enumi})}}
\item $F^\beta + Z^\beta: \ {\cal H}_1 \oplus {\cal H}_2 \ \to \
{\cal H}_1 \oplus {\cal H}_2$ is injective,
\item $Z^\beta$ is of finite rank and hence compact,
\item $Z^\beta \ = \ - G^\beta Z^\beta$,
\item for all $\gamma \in [\beta_1,\beta_2]$ we have ${\rm im}\
Z^\beta \subset
{\cal D}(E^\gamma)$, and $E^\gamma Z^\beta$ is compact and independent
of $\gamma$, i.e.
\[ E^\gamma Z^\beta \ = \ E^\beta Z^\beta.\]
\end{enumerate}
\end{lemma}
\proof For (i), (ii), (iii) cf. \cite{B2}. For (iv), we first prove
${\rm im}\ G^\beta \subset {\cal D}(E^\gamma)$. $f \in {\rm im}\ (R_i+\psi_2'
P_{{\rm bd},\beta} \psi)$ implies $f \in {\cal D}\bigl( (D^\gamma)^* \bigr)$
by using the following facts:
\begin{itemize}
\item $R_i$ is infinitely smoothing,
\item $\psi_2' P_\beta \psi$ maps ${\cal H}_2$ into ${\cal
D}(D^\beta)$,
\item $f \in {\rm im}\ (R_i+\psi_2' P_{{\rm bd}, \beta} \psi)$ implies $f \vert
(0,\epsilon_1) = 0$,
\item $D^\beta$ and $(D^\gamma)^*$ are both closed extensions of
first order elliptic differential operators.
\end{itemize}
Analogously, we obtain ${\rm im}\ ({L_i}^*+\psi' {P_{{\rm bd},\beta}}^* 
\psi_2)
\subset {\cal D}(D^\gamma)$. Thus we conclude with (iii):
\[ {\rm im} \ Z^\beta \ \subset \ {\rm im}\ G^\beta \ \subset \ {\cal
D} (E^\gamma). \]
Obviously, the considerations above imply also $E^\gamma Z^\beta =
E^\beta Z^\beta$. $E^\beta Z^\beta$ is compact, since the restriction of
$E^\beta$ to the finite dimensional vector space ${\rm im}\ Z^\beta$
is bounded. \hfill\endproof

\noindent
Now, $E^\gamma(F^\gamma + Z^\gamma) = I + F^\gamma + G^\gamma + E^\beta
Z^\beta$ is a Fredholm operator with index $0$. This implies together
with Lemma \ref{propz} (i), that $E^\beta(F^\beta + Z^\beta)$ is
bijective and hence invertible. Lemma \ref{pcont} guarantees the
continuity of
\[ [\beta_1,\beta_2] \ni \gamma \ \mapsto \ E^\gamma(F^\gamma +
Z^\gamma) \ =
\ I + F^\gamma + G^\gamma + E^\beta Z^\beta, \]
hence $\bigl( E^\gamma ( F^\gamma + Z^\beta) \bigr)^{-1}$ is well
defined for $\gamma$ close to $\beta$ and continuous in $\gamma$. This
finally implies continuity of the map
\[ \gamma \ \mapsto \ (E^\gamma)^{-1} \ = \ (F^\gamma + Z^\beta) \bigl(
E^\gamma (F^\gamma + Z^\beta) \bigr)^{-1} \]
and finishes the proof of ${\rm ind}\ D^{\beta_1} \ = \ {\rm ind}\
D^{\beta_2}$. Let us state this result in a theorem:
\begin{theorem} \label{equalind}
Let $1 \le \beta_1 < \beta_2$ and $\{ D^\beta \}_{\beta \in
[\beta_1,\beta_2]}$ a family of {\em unperturbed operators}
which are of the type {\myref{operators}} and coincide on
$M_\epsilon$. Moreover, we assume $h_\beta$ to be strictly
increasing, $h_\beta\ \vert\ (0,\epsilon_0) =
x^\beta$ and $h_\beta$ converge uniformly on
compact intervals for $\beta \to \gamma$. Then
\[ {\rm ind}\ D^{\beta_1}_\delta \ = \ {\rm ind}\ D^{\beta_2}_\delta. \]
\end{theorem}
In the sequel we will apply the results of this and the previous
section to geometric operators. In the case of the Spin--Dirac
operator we will prove that the index does not change by deforming
metric horns into metric cones. However, in the case of the 
Gau{\ss}-Bonnet or the Signature operator we may only prove
constance of the index for
deformations from horns with warping function $x^{\beta_2}$ to horns
with warping function $x^{\beta_1}$ (both $\beta_1, \beta_2 > 1$). This
is due to the fact, that
there occurs a ``perturbation'' in the latter two operators which can
only be removed for $\beta > 1$ (cf. Proposition
\ref{removeperturbation}).
Fortunately, we can derive a Gau{\ss}-Bonnet theorem for metric
horns by $L^2$-cohomology arguments of J. Cheeger.


\section{Application to the Spin--Dirac operator}
Let $M$ be an open spin manifold of even dimension $m$ and let
\[ h(x),\ D,\ D^+,\ S(M)\ {\rm and}\ S(N) \]
be as described in Section 2 in the example of the Spin--Dirac operator.
We are interested in
a classification of all closed extensions $\bar D^+$ of
$D^+$ and
a formula for ${\rm ind}\ \bar D^+$ in the spirit of \cite{APS}.
The answers to these questions in the regular singular case,
i.e. $h(x) = x$ on $(0,\epsilon_0)$ are due to A.W. Chou, which we recall
for the sake of completeness. Here we use the terminology of \cite{BS}:
\begin{theorem} {\rm (see \cite[Theorems (3.2) and (5.23)]{Chou})}
\label{chou} Let $h(x)=x$ on $(0,\epsilon_0)$.
The closed ex\-ten\-sions of $D^+$ are in
$1:1$ correspondence to subspaces $W$ of
\[ \bigoplus_{\vert s \vert < {1 \over 2}}\ {\rm ker}\ (D_N - s) \]
and
\begin{equation}
{\rm ind}\ D^+_W \ = \ \int_M \hat A \ - \ {1 \over 2} \bigl(
\eta(0) + b \bigr) \ + \ \underbrace{ {\rm dim}\ W \ - \ \sum_{-{1\over
2} < s < 0} {\rm dim}\ E_s}_{(*)},
\end{equation}
where $\hat A$ is the Hirzebruch $\hat A${\em -polynomial},
$\eta(0)$ the eta-invariant of $D_N$,
$E_s$ the eigenspace of $D_N$ to the
eigenvalue $s$ and $b:={\rm dim\ ker}\ D_N$.
For the particular extension $D^+_\delta$ as
defined in \cite{BS}, the term $(*)$ vanishes.
\end{theorem}
In Theorem \ref{chou}, we chose
the orientation of \cite{APS} which is opposite to the orientation of
\cite{Chou}. This causes different signs of $\eta(0)$ in the index
formulas.

The index formula of Chou can be rederived from the index
formula in
Theorem 4.1. of \cite{BS} by using the identification \myref{sdo}. The
vanishing of the residua in (4.52) of \cite{BS} is guaranteed by
\cite[Theorem (4.2)]{APS}.

Now we turn over to the case of metric horns and assume $h(x) = x^\alpha$ on
$(0,\epsilon_0)$ with $\alpha > 1$. The idea is to establish a
connection
between horns and cones and to transfer the known results for cones to
horns. The following technical lemma is useful to set up this connection:

\begin{lemma} \label{hexist}
Let $0 < \epsilon \le 1$, $1 \le \beta_1 < \beta_2$ and $h_{\beta_1},
h_{\beta_2} \in C^\infty((0,\epsilon))$ with the following properties:

\medskip\noindent
\widearray
\begin{tabular}{rll}
\hphantom{p}--& $h_{\beta_j}(x)\ =\ x^{\beta_j}$ &on
$(0,\epsilon_0)$,\\
--& $h_{\beta_j}'(x)\ \ge\ c \ >\ 0$ &on $[\epsilon_0,\epsilon_1]$,\\
--& $h_{\beta_1}(x) \ = \ h_{\beta_2}(x)$ & on
$(\epsilon_1,\epsilon)$.\\
\end{tabular}
\medskip

\noindent
Then there exists a family $\{ h_\beta \}_{\beta \in [\beta_1,\beta_2]}
\subset C^\infty(0,\epsilon)$ and an $\tilde \epsilon \in
(0,\epsilon_0)$ such that

\medskip\noindent
\widearray
\begin{tabular}{rll}
\hphantom{p}--&
    $h_\beta(x)\ =\ x^\beta$& on $(0,\tilde\epsilon)$,\\
--& $h_{\beta}'(x)\ \ge\ c_1 \ >\ 0$ & on $[\tilde
\epsilon, \epsilon_1]$,\\
--& $h_\beta(x)$ independent of $\beta$ & on
$(\epsilon_1,\epsilon)$,\\
--& $h_\beta \Rightarrow h_\gamma$ and $h_\beta' \Rightarrow
h_\gamma'$ & on compact intervals as $\beta \to \gamma$.
\end{tabular}
\end{lemma}
\proof $h_\beta$ is given by
\[ h_\beta(x) \ := \ \varphi(x) x^\beta \ + \ (1-\varphi(x)) \bigl( a_\beta
h_{\beta_1}(x) + (1-a_\beta) h_{\beta_2}(x) \bigr), \]
where $\varphi \in C^\infty(0,\epsilon)$ is monotonic decreasing,
$\varphi \vert (0,\tilde \epsilon) \equiv 1$, $\varphi \vert
(\epsilon_0,\epsilon) \equiv 0$ and
$\displaystyle a_\beta := {{\epsilon_0}^\beta - {\epsilon_0}^{\beta_2}
\over {\epsilon_0}^{\beta_1} - {\epsilon_0}^{\beta_2}}$.\endproof

Let $(M,g^\alpha)$ be a spin manifold with a metric horn. The previous lemma
enables us to establish a metric homotopy $(M,g^\beta)$ from a horn
$(\beta = \alpha)$ to a cone $(\beta = 1)$ by assuming
\[ g^\beta \vert U \ \cong\ dx^2 \ + \ {h_\beta}^2\ g_N. \]
Obviously all conditions are fulfilled to apply Theorem \ref{equalind}
with $(\beta_1,\beta_2) = (1,\alpha)$.
Hence we conclude from Chou's index formula:
\[ {\rm ind}\ D^{+,\alpha}_\delta \ = \ {\rm ind}\ D^{+,1}_\delta \
= \ \int_{M_\epsilon} \hat A \ + \
\int_{(U,g^1)} \hat A \ - \ {1 \over 2} (\eta (0) + b). \]
The same holds true for manifolds with more than one metric horn and
different warping exponents $\alpha_j > 1$ for each horn. A
drawback of this index formula yet is the occurence of the term
$\int_{(U,g^1)} \hat A$. Fortunately, the $O(m)$-invariance of $\hat A$
implies $\hat A_{m/2} \equiv 0$ on any warped product $dx^2+ h(x)^2 g_N$.
Here, $\hat A_{m/2}$ denotes the homogeneous component of $\hat A$ of degree
$m/2$. Using these facts we conclude from the last two sections 
(note, that for the
Spin--Dirac operator $H_{1 \over 2\alpha} = H = {\rm ker}\ D_N$):
\begin{theorem} \label{spindirac}
Let $M$ be a singular spin manifold with metric horns.
The closed extensions of the Spin--Dirac operator
$D^+$ are in $1:1$ correspondence to subspaces $W \subset
{\rm ker}\ (D_N)$ and
\[ {\rm ind}\ D^+_W \ = \ \int_M \hat A \ - \ {\eta(0) \over 2} \ + \
\bigl( {\rm dim}\ W - {b \over 2} \bigr). \]
In particular, for the minimal and maximal extension we have
\[ {\rm ind}\ D^+_{\rm max/min} \ = \ \int_M \hat A \ - \ {\eta(0)
\over 2} \ \pm \ {b \over 2}.\]
\end{theorem}

An easy consequence of the
Theorem of Lichnerowicz (cf. \cite[Corollary II.8.9.]{LM}) is
\begin{cor} Let $M$ be as in Theorem {\rm \ref{spindirac}} and the scalar
curvature of $g_N$ be positive everywhere. Then
the Spin--Dirac operator on $M$ has a unique closed extension.
\end{cor}
These results allow the conclusion, that the case of
metric horns is considerably easier than the case of metric cones:
No longer all eigenvalues in $(-{1 \over 2},{1 \over 2})$ of $D_N$ are of
importance but only its kernel. The latter is independent of
conformal changes of the metric $g_N$ (see \cite{H}) whereas the small
eigenvalues are not stable under even such an easy deformation as
multiplying the metric $g_N$ by a positive constant.


\section{Gau{\ss}-Bonnet and Signature operator}
According to Section 3, the closed extensions of
the Gau{\ss}-Bonnet and the Signature operator on manifolds with metric
horns are characterized by subspaces of
$H_{1 \over 2\alpha} = \oplus_{\vert \alpha s \vert < 1/2}
{\rm ker}\ (S-s)$. The examples in Section 2 describe the concrete choices
of $S$, $H$, $\tilde S$, $\tilde H$ and $\tilde A$ for both operators.
In either case we choose $H$ equals the space of
harmonic forms ${\cal H}(N)$ and $S$ to be a diagonal matrix operator with
diagonal entries given by $c_j := (-1)^j (j - {n \over 2})$ or
$b_j := j - {n \over 2}$, respectively. From this we conclude easily
\begin{theorem} \label{clext}
Let $M$ be an oriented singular Riemannian manifold with metric horns and
$D_{\rm GB}$ and $D_{\rm S}$ the Gau{\ss}-Bonnet operator and the Signature
operator on $M$ (for the latter we assume ${\rm dim}\ M = 4k$).
\renewcommand{\labelenumi}{{\rm \arabic{enumi}.}}
\begin{enumerate}
\item If ${\rm dim}\ M$ is even, $D_{\rm GB}$ has a unique closed
extension.
\item If ${\rm dim}\ M$ is odd, the closed extensions of $D_{\rm GB}$
can be characterized by the subspaces of ${\cal H}^{n/2}(N)$, where $N$ is
the cross section of the horns and $n = {\rm dim}\ N$.
\item $D_{\rm S}$ has a unique closed extension.
\end{enumerate}
Moreover, all closed extensions are Fredholm.
\end{theorem}
Now, let us consider the indices of $D_{\rm GB}$ and
$D_{\rm S}$ on metric horns. Index formulas for the Gau{\ss}-Bonnet
and the Signature operator on manifolds with metric cones were first
obtained by J. Cheeger (see \cite{Che2}). Similar formulas for the 
asymptotic conic case are derived in \cite{BS} (see Theorem 5.1. and
5.2.).
>From the results in Section 4 we may deduce
\begin{theorem} \label{indgbands}
Let $M$ be an even dimensional oriented open manifold. Let $g_1$ and $g_2$
be two metrics on $M$ which induce the structure of a singular
Riemannian manifold with a metric horn for the same choice of
$U$ and identification $U \cong (0,\epsilon) \times N$.
Moreover, we assume $g_j \vert U \cong dx^2 + {h_{\beta_j}}^2 g_N$ and $g_1
\vert M_{\epsilon_1} = g_2 \vert M_{\epsilon_1}$ for a suitable
$\epsilon_1 \in
(0,\epsilon)$ and the following properties of the warping functions
$h_{\beta_j}$:
\begin{itemize}
\item $h_{\beta_j}(x) \ = \ x^{\beta_j}$ near the singularity and $\beta_1,
\beta_2 > 1$,
\item $h_{\beta_j}$ is strictly increasing on $(0,\epsilon_1)$.
\end{itemize}
Then the unique closed extensions of the Gau{\ss}-Bonnet
resp. Signature operators $D_{\rm GB/S}^{\beta_j}$ have the same
index:
\[ {\rm ind}\ D_{\rm GB/S}^{\beta_1} \ = \ {\rm ind}\ D_{\rm
GB/S}^{\beta_2}. \]
Of course, for the Signature operator we have additionally to assume
${\rm dim}\ M = 4k$.
\end{theorem}
This theorem generalizes to manifolds
with several horns
in an obvious manner. However, in contrast to the Spin--Dirac operator,
our method does not allow
to conclude coincidence of the indices for metric horns
and metric cones. This lack is due to the occurence of a perturbation
$\tilde A \not\equiv 0$ in these two operators.
\proof As in the last section, Lemma \ref{hexist} implies the existence of
metrics $\{ g^\beta \}_{\beta \in [\beta_1,\beta_2]}$ on $M$ being a homotopy
connecting $g_1$ with $g_2$. By Proposition \ref{removeperturbation} we
conclude for each $\beta \in [\beta_1,\beta_2]$:
\begin{equation} \label{eqdomain}
{\cal D}(D_{{\rm GB},0}^\beta) \ = \ {\cal D}(D_{\rm GB}^\beta) \quad
{\rm and}
\quad {\rm ind}\ D_{{\rm GB},0}^\beta\ =\ {\rm ind}\ D_{\rm GB}^\beta
\end{equation}
where $D_{{\rm GB},0}^\beta$ is the ``unperturbed version'' of $D_{\rm
GB}^\beta$. For the unperturbed operators we conclude with Theorem
\ref{equalind}:
\begin{equation} \label{ggg}
{\rm ind}\ D_{{\rm GB},0}^{\beta_1}\ =\
{\rm ind}\ D_{{\rm GB},0}^{\beta_2}.
\end{equation}
\myref{eqdomain} and \myref{ggg} together imply the statement of the 
theorem for the Gau{\ss}-Bonnet operator and analogously for the Signature
operator. Note that in contrast to \myref{eqdomain} the domains
of the two operators in \myref{ggg} generally do not coincide. \endproof

In the case of the Gau{\ss}-Bonnet operator one may easily 
conclude a stronger result from $L^2$-cohomology 
considerations of J. Cheeger.
These $L^2$-cohomology considerations imply for example, that
the statement of Theorem \ref{indgbands}
holds also for $\beta_1 = 1$.
Let us explain this in more detail: Let $M$, $g_1$ and $g_2$ be as in
Theorem \ref{indgbands} with the only difference $\beta_1 = 1$.
Using the terminology
of {\em Hilbert complexes} (see \cite{BL1} 
for a detailed treatment of this notion) we may conclude for both metrics, 
that each $d_j: \Omega_0^j(M) \to \Omega_0^{j+1}(M)$ has a unique closed
extension $D_j: {\cal D}(D_j) \to L^2(\Lambda^{j+1}T^*M)$ and that
\[0 \ \to \ {\cal D}(D_0) \ \stackrel{D_0}{\to} \ {\cal D}(D_1) \
\stackrel{D_1}{\to} \ \cdots \ {\cal D}(D_{m-1}) \ \stackrel{D_{m-1}}{\to} \
L^2(\Lambda^m T^*M) \ \to \ 0\]
is a Fredholm complex (cf. e.g. \cite[Theorem 3.7.(a)]{BL2}
for the metric $g_1$).
In the case of the metric $g_2$ this follows
easily from uniqueness and the Fredholm property of the closed
extension of $D_{\rm GB}^{\beta_2}$. Theorem 2.1. of \cite{Che1} describes
the $L^2$-cohomology groups  ${\rm ker}\ D_j / {\rm im}\ D_{j-1}$
in terms of relative cohomology groups of $(M_\epsilon,N)$. 
Being stated for manifolds with metric cones, Theorem 2.1.
also holds true for horns since the only necessary
tools ($L^2$-versions of 
Poincar{\'e} lemma and Mayer Vietoris developed in \cite{Che0}) are valid for
both situations. This implies that corresponding $L^2$-Betti numbers are
exactly the same for both metrics $g_1$ and $g_2$. In particular,
the $L^2$-Euler characteristics coincide.
Denoting by $D_{\rm GB}^1$ and $D_{\rm GB}^{\beta_2}$ the Gau{\ss}-Bonnet 
operators corresponding to the metrics $g_1$ and $g_2$, we conclude with
\cite[Theorem 2.4]{BL1} and \cite[Theorem 3.7(d)]{BL2}:
\begin{equation} \label{indagree}
{\rm ind}\ D_{\rm GB}^{\beta_2} \ = \ \left\{ \begin{array}{cc} {\rm ind}\ D_{\rm
GB, max}^1 & {\rm if}\ {m \over 2}\ {\rm is\ even}, \\ & \\
{\rm ind}\ D_{\rm GB, min}^1 & {\rm if}\ {m \over 2}\ {\rm is\ odd}. 
\end{array} \right.
\end{equation}
Note, that the distinction 
between $D_{\rm GB,min}$ and $D_{\rm GB, max}$ for the conic
metric $g_1$ is necessary: though there is a unique ideal boundary
condition there may be many closed extensions of
the Gau{\ss}-Bonnet operator 
$D_{\rm GB}: \Omega_0^{\rm even} \to \Omega_0^{\rm odd}$. 

\myref{indagree} establishes a connection between metric horns and cones and
allows to proceed in the same way as in the previous section. Henceforth we
change our notations from $\beta_2$ to $\alpha$ and from $g^2$ to 
$g^\alpha$. Moreover, choose $\epsilon_0 \in (0,\epsilon)$ such that
$h_1$ resp. $h_\alpha$ coincide with $x$ resp. $x^\alpha$ on $(0,\epsilon_0)$.
With the help of the Gau{\ss}-Bonnet formula for metric cones (see
\cite[{[6.1]}]{Che1} or \cite[Theorem 5.1]{BS}) we conclude:


\begin{eqnarray}
  \DST {\rm ind}\ D_{\rm GB}^\alpha &=&\DST
 \int_{(M_{\epsilon_0},g^1)} e \ + \ \int_{(U_{\epsilon_0},g^1)} e\nonumber \\[1.2em]
  &&\DST+ \ {1 \over 2} \Bigl( \sum_{j=0}^{{m \over 2}-1} (-1)^j b_j(N) \ - \ 
 \sum_{j={m \over 2}}^n (-1)^j b_j(N) \Bigr)  \label{preindgb}\\[1.2em]
  &&\DST+\ \underbrace{ \sum_{p \ge 1} \alpha_p\ {\rm Res}_1\ 
 \eta_{S_0}(2p) \ - \ {\eta_3(0) \over 2} }_{=: \ \Xi(N)}. \nonumber
\end{eqnarray}

Here $e$ denotes the Euler-class and $b_j(N)$ the Betti numbers of the cross
section $N$. The
ingredients $\eta_{S_0}$ and $\eta_3(0)$ are defined as in \cite{BS}. We
emphasize that $\Xi(N)$ is made up of spectral data on the cross section 
$N$, which generally are difficult to calculate.  

The following two lemmas allow us to express the right-hand side of
\myref{preindgb} completely in terms of the Riemannian manifold
$(M,g^\alpha)$. 
\begin{lemma} \label{gbforcollars}
Let $(N,g_N)$ be an odd dimensional compact manifold, $n = {\rm dim}\ N$
and $U_{\delta \epsilon}$ the warped
collar $[\delta,\epsilon] \times N$ with metric $g_{U_{\delta \epsilon}} 
= dx^2 + h(x)^2 g_N$.
Then there exist differential forms 
$\alpha_k$ ($k = 1,2,\dots,{n-1 \over 2}$) on $N$ depending
only on the intrinsic metric $g_N$, such that
\[ \int_{U_{\delta \epsilon}} e \ = \ 
\sum_{k=0}^{n-1 \over 2} \Bigl( (h'(\epsilon))^{2k+1} - 
(h'(\delta))^{2k+1} \Bigr) \int_N \alpha_k. \]
Consequently, $\int_{U_{\delta \epsilon}} e = 0$ for linear warping
functions $h$. 
\end{lemma}
\proof The
Chern-Gau{\ss}-Bonnet theorem for manifolds with boundary yields
\[ \int_{U_{\delta \epsilon}} e \ = \ 
\underbrace{\chi\Bigl( [\delta,\epsilon] \times N \Bigr)}_{= 0} \ -
\ \int_{\partial U_{\delta \epsilon}} Se \]
where $Se$ is a $SO$-invariant form defined near the boundary of 
$U_{\delta \epsilon}$. According to \cite[pp. 252]{G}, 
$Se$ can be written locally with 
respect to an oriented 
orthonormal frame $\{ v_j \}$ of $TU_{\delta \epsilon}$, 
$v_{n+1} = {\partial \over \partial x}$ as 
\begin{equation} \label{Se}
Se \ = \ \sum_{k=0}^{n-1 \over 2} c_{k,m} \sum_{\sigma \in {\cal S}_n}
{\rm sgn}(\delta) \Omega^U_{\sigma_1 \sigma_2} \wedge \dots \wedge
\Omega^U_{\sigma_{2k-1} \sigma_{2k}} \wedge \omega^U_{\sigma_{2k+1} n+1}
\wedge \dots \wedge \omega^U_{\sigma_n n+1}.
\end{equation}
$c_{k,m}$ are suitable chosen constants. Let $\omega^N$, $\omega^U$ and
$\Omega^N$, $\Omega^U$ be connection and curvature forms w.r.t. 
orthonormal frames $e_1,\dots,e_n$ and ${1 \over h}e_1, \dots,
{1 \over h}e_n, {\partial \over \partial x}$ of $N$ and $U$. A
standard calculation yields: 
\begin{eqnarray*}
\omega^U_{j k} &=& \left\{ \begin{array}{cc} \omega^N_{j k} & {\rm for}\ 1
\le j,k \le n, \hfill \\ & \\
- h'\ e^j & {\rm for}\ k=n+1 \ {\rm and}\ 1 \le j \le n, \end{array} \right. \\
\Omega^U_{j k} &=& \Omega^N_{j k} \ - \ (h')^2\ e^j \wedge e^k \qquad {\rm
for}\ 1 \le j,k \le n.
\end{eqnarray*}
Inserting these identities in \myref{Se} one easily deduces how to choose the
intrinsic defined forms $\alpha_k$ on $N$. \endproof
\begin{lemma}[{\rm \cite[p. 607]{Che2}}] \label{Chelemma}
Using the notations above, we have
\begin{equation} \label{Chel1}
\Xi(N) \ = \ \int_N\ \sum_{k=0}^{n-1 \over 2} 
\bigl( h_1'(\epsilon_0) \bigr)^{2k-1} \alpha_k.
\end{equation}
\end{lemma}
Cheeger proved this identity between the spectral invariant $\Xi(N)$
and the boundary integral by comparing his
Gau{\ss}-Bonnet formula for metric cones 
with the classical Gau{\ss}-Bonnet
formula for manifolds with boundary. Using the notions of \cite{BS},
$\Xi(N)$ has the form
\begin{eqnarray*} 
\Xi(N) &=& \underbrace{\sum_p 
\Bigl( \alpha_p\ {\rm Res}_1\ \eta_{S_0}(2p) \ - \ 
{\beta_p \over 2} \sum_k (-1)^k\ {\rm Res}_1\ \zeta_k(2p+1) \Bigr)}_{(*)} \\
&+& {1 \over 2} \sum_k (-1)^{k+1}\ {\rm Res}_1\ \zeta_k(1).
\end{eqnarray*}
Actually, Cheeger states the identity \myref{Chel1} without the term $(*)$ in
$\Xi(N)$, which is most likely equal $0$. 
However, since our considerations are based on \cite{BS} we
include this term to establish consistency. Let us sketch Cheeger's proof.
\proof $L^2$-cohomology arguments imply
\begin{equation} \label{coho1}
\chi_{(2)}(M) \ = \ \chi(M_\epsilon) \ + \ {1\over 2}\chi_{(2)}(CN) \ + \
{1\over 2}\chi_{(2)}(CN,N),
\end{equation}
where $CN$ denotes the metric cone $(0,\epsilon_0) \times N$, and
\begin{equation} \label{coho2}
\chi_{(2)}(CN) \ = \ \sum_{j=0}^{{m\over 2}-1} (-1)^j b_j(N), \qquad
\chi_{(2)}(CN,N) \ = \ -\ \sum_{j={m\over 2}}^n (-1)^j b_j(N).
\end{equation}
Since ${\rm dim}\ M_\epsilon$ is even we do not have to bother about
absolute or relative boundary conditions for the corresponding 
Euler characteristics are the same (see \cite[Theorem 4.2.7]{G}).
Using \myref{coho1}, \myref{coho2} together with 
the Gau{\ss}-Bonnet formula for cones in the terminology of \cite{BS}
we conclude
\[ \chi(M_{\epsilon_0}) \ = \ \int_{(M_{\epsilon_0},g^1)} e \ + \ \Xi(N). \]
On the other hand the classical Chern-Gau{\ss}-Bonnet formula for
$M_\epsilon$ reads as
\[ \chi(M_\epsilon) \ = \ \int_{(M_{\epsilon_0},g^1)} e \ + \ \int_{\partial
M_{\epsilon_0}} Se. \]
Similarly to the previous proof we deduce
\[ \int_{\partial M_{\epsilon_0}} Se \ = \ \int_N\ \sum_{k=0}^{n-1 \over 2}
\bigl(h_1'(\epsilon)\bigr)^{2k+1} \alpha_k, \]
which finishes the proof. \endproof 

With all these identities at hand we conclude
\begin{theorem} {\rm (Gau{\ss}-Bonnet formula for metric horns)}
\label{gaussbonnet} Let $M$ be a singular manifold with metric horns 
of even dimension $m=n+1$. There exists a unique closed extension of
$D_{\rm GB}: \Omega^{\rm even}_0(M) \to \Omega^{\rm odd}_0(M)$
which is Fredholm and its index is given by
\[ {\rm ind}\ D_{\rm GB}\ = \
\int_M e \ + \ {1 \over 2}\Bigl( \sum_{j=0}^{{m\over 2}-1}
(-1)^j b_j(N) \ - \ \sum_{j={m\over 2}}^n (-1)^j b_j(N) \Bigr). \]
\end{theorem}
\proof Let $\delta\in (0,\epsilon_0)$ be arbitrary.
Using Lemma \ref{gbforcollars} twice and the fact that $h_\alpha$ and
$h_1$ agree near $x=\epsilon$ we obtain:
\begin{equation} \label{help1}
\int_{(U_{\delta \epsilon},g^1)} e \ = \ \int_{(U_{\delta
\epsilon},g^\alpha)} e \ + \ \int_N \ \sum_k \Bigl( \bigl( h'_\alpha(\delta)
\bigr)^{2k+1} \ - \ \bigl( h'_1(\delta) \bigr)^{2k+1} \Bigr)\ \alpha_k.
\end{equation}
Since $h_1$ is linear on $(0,\epsilon_0)$: $\displaystyle
\int_{(U_{\delta},g^1)} e = 0$.
Together with Lemma \ref{Chelemma}, \myref{preindgb},
and \myref{help1}, this implies
\begin{eqnarray*}
{\rm ind}\ D_{\rm GB}^\alpha &=&
\underbrace{\int_{(M_{\delta},g^\alpha)} e \ + \int_N \ \Bigl( \sum_k
\bigl(h_\alpha'(\delta)\bigr)^{2k+1} \ \alpha_k \Bigr)}_{(*)} \\
&+& {1 \over 2} \Bigl( \sum_{j=0}^{{m \over 2}-1} (-1)^j b_j(N) \ - \ 
\sum_{j={m\over 2}}^n  (-1)^j b_j(N) \Bigr).
\end{eqnarray*}
Using the fact that $h_\alpha'(\delta)
\to 0$ as $\delta \to 0$ we see that $(*)$ coincides with $\displaystyle
\int_{(M,g^\alpha)} e$. \endproof 

Let us consider the following situation as an application of the previous 
results: Let $(M,g^1)$ be an oriented, compact closed manifold of even
dimension $m$ and $p_1, \dots, p_k$ arbitrary points of $M$. Without loss of
generality we may assume, that in the neighborhood of each point $p_j$
the metric $g^1$ is of the form $dx^2 + x^2 g_{N_j}$ with $N_j$ diffeomorphic
to the sphere $S^{m-1}$. We may consider $M_0 := M-\{ p_0, \dots, p_k\}$ 
as well as a singular manifold with metric cones. 
Now, let the metric change continuously so that the
neighborhoods of each $p_j$ become  metric horns. We denote the family of
metrics on $M_0$ again by $\{ g^\beta \}_{\beta \in [1,\alpha]}$. Whereas
the index of the Gau{\ss}-Bonnet operator (corresponding to the unique ideal
boundary condition) does not change, the map  
\[ [1,\alpha] \ni \beta \ \to \ \int_{(M_0,g^\beta)} e \in {\Bbb{Z}} \]
skips from $\chi(M)$ to $\chi(M)-k$ as soon as $\beta$ becomes
greater than $1$. This is somewhat surprising since the metric of $M_0$
changes smoothly. Such a phenomenon does not occur for the $O(m)$-invariant
forms $\hat A$ or $L_k$ (Hirzebruch's L-polynomial with $k := m/4$).
This is an easy consequence of the fact 
that the integral of $O(m)$-invariant forms over warped products vanishes as
shown in the last section. 




Unfortunately, our method only partially applies to
the Signature operator. By Theorem \ref{clext} the Signature
operator has a unique closed extension which is Fredholm.
However, we are not able to prove an index theorem for the
Signature operator. By analogy the following result is conceivable.
It was stated as a conjecture in an earlier version of this paper.
Now it is a Theorem since a proof has been 
announced by J. Br\"uning \cite{B3}. 

\begin{theorem} {\rm (Signature formula for metric horns)}
\label{signature} Let $M$ be a $4k$-di\-mensional singular manifold
with metric horns. There exists a unique closed extension of
$D_{\rm S}: \Omega^+_0(M) \to \Omega^-_0(M)$ which is Fredholm and its index
is given by
\[ {\rm ind}\ D_{\rm S}\ = \ \int_M L_k \ - \ \eta(N), \]
where $L_k$ is the $k$-th Hirzebruch $L$-polynomial
and $\eta(N)$ is the eta-invariant of the operator
\[ \alpha \ \mapsto \ (-1)^{k+j+1}(*_N d_N - d_N *_N) \quad {\rm for}\
\alpha \in \Omega^{2j}(N), \]
and $N$ is the cross section of the horns.
\end{theorem}

Br\"uning's method is different from ours. He announces a heat trace
asymptotics for metric horns and hence he does not reduce the problem
to the conic case. However, he needs our Theorem \ref{clext}.


It may be useful to note, that the Signature operator on manifolds with horns
decomposes into the infinite direct sum of operators acting on one and two
dimensional subspaces. 
It was our hope that this could be used to prove the graph continuity
(cf. Theorem \ref{equalind}) for $\beta\to 1$ for the Signature operator,
too. Unfortunately, we did not succeed.
However, this decomposition seems to be of some interest in its own.       
An analogous decomposition of the Gau{\ss}-Bonnet
operator was proven in \cite[Lemma 2.2.]{BL2}.

We decompose $L^2(\Lambda^* T^*N)$ into
\begin{equation} \label{dec1}
L^2(\Lambda^* T^*N)\ = \ {\cal H} \oplus H_1 \oplus H_2 \oplus H_3
\oplus H_4 \oplus H_5 \oplus H_6,
\end{equation}
where ${\cal H}$ denotes the space of harmonic forms on $N$ and $H_l$
are chosen as follows (the matrices $M_l$ will be used
below):

\medskip

\noindent
\begin{center}
\begin{tabular}{|c|c|c|} \hline
\rule[-1mm]{0mm}{7mm} $l$ & $H_l$ & $M_l$ \\ \hline \hline
\rule[-5mm]{0mm}{12mm} $1$ & $\displaystyle \bigoplus_{0 \le j \le k-2
\atop \lambda>0}
\Bigl( E^{2k-2j-2}_{\lambda,cl} \oplus E^{2k+2j+2}_{\lambda,cl}
\Bigr)$ & $\displaystyle \left( \begin{array}{cc} (2j+{3 \over 2})h' &
(-1)^{j+1}\sqrt{\lambda} \\ (-1)^{j+1}\sqrt{\lambda} &
-(2j+{5 \over 2})h' \end{array} \right)$ \\ \hline
\rule[-5mm]{0mm}{12mm} $2$ & $\displaystyle \bigoplus_{0 \le j \le k-2 \atop
\lambda > 0} \Bigl( E^{2k-2j-3}_{\lambda,ccl} \oplus E^{2k+2j+1}_{\lambda,ccl}
 \Bigr)$ & $\displaystyle \left( \begin{array}{cc} (2j+{5 \over 2})h' &
(-1)^{j+1}\sqrt{\lambda} \\ (-1)^{j+1}\sqrt{\lambda} &
-(2j+{3 \over 2})h' \end{array} \right)$ \\ \hline
\rule[-5mm]{0mm}{12mm} $3$ & $\displaystyle \bigoplus_{0 \le j \le k-1
\atop \lambda>0} \Bigl(
E^{2k-2j-1}_{\lambda,cl} \oplus E^{2k+2j+1}_{\lambda,cl} \Bigr)$ &
$\displaystyle \left( \begin{array}{cc} (2j+{1 \over 2}h' &
(-1)^j\sqrt{\lambda} \\ (-1)^j\sqrt{\lambda} & -(2j+{3 \over 2}h' \end{array}
\right)$ \\ \hline
\rule[-5mm]{0mm}{12mm} $4$ & $\displaystyle \bigoplus_{0 \le j \le k-1
\atop \lambda>0} \Bigl(
E^{2k-2j-2}_{\lambda,ccl} \oplus E^{2k+2j}_{\lambda,ccl} \Bigr)$ &
$\displaystyle \left( \begin{array}{cc} (2j+{3 \over 2}h' &
(-1)^j\sqrt{\lambda} \\ (-1)^j\sqrt{\lambda} & -(2j+{1 \over 2})h'
\end{array} \right)$ \\ \hline
\rule[-5mm]{0mm}{12mm} $5$ & $\displaystyle \bigoplus_{\lambda > 0}
E^{2k}_{\lambda,cl}$ & {\ } \\ \hline
\rule[-5mm]{0mm}{12mm} $6$ & $\displaystyle \bigoplus_{\lambda > 0}
E^{2k-1}_{\lambda,ccl}$ & {\ } \\ \hline
\end{tabular}
\end{center}
\medskip

$E_{\lambda,cl/ccl}^j$ denotes the space of closed resp. coclosed
$j$-eigenforms of the Laplacian on $N$ corresponding to the eigenvalue 
$\lambda$. 

$E^{2k}_{\lambda,cl}$ resp. $E^{2k-1}_{\lambda,ccl}$ 
admit a further decomposition (into $\pm 1$-eigenspaces) via the
involutions ${1 \over \sqrt{\lambda}}d_N *_N$ resp. ${1 \over
\sqrt{\lambda}} *_N d_N$:
\begin{equation} \label{dec2}
E^{2k}_{\lambda,cl} \ = \ E^{2k,+}_{\lambda,cl} \ \oplus \
E^{2k,-}_{\lambda,cl} \qquad {\rm and} \ \ E^{2k-1}_{\lambda,ccl} \ = \
E^{2k-1,+}_{\lambda,ccl} \ \oplus \ E^{2k-1,-}_{\lambda,ccl}.
\end{equation}
By arguments analogously to \cite{BL2} we conclude
\begin{lemma} \label{sepds} The operator $D_{\rm S}$ reduces with respect
to the decompositions {\rm \myref{dec1}, \myref{dec2}} in the neighborhood of a
metric horn into
\begin{eqnarray*}
D_{\rm S} \vert H &\cong& \partial_x + {h' \over h} S_0 \ = \
\bigoplus_{j=0}^n \bigl( \bigoplus_{{\rm dim\ } {\cal H}^j} ( \partial_x + {1
\over h} b_j ) \bigr),\\
D_{\rm S} \vert V &\cong& \partial_x + {1 \over h} M_l
\quad {\rm w.r.t.}\ \{ \eta_1, \eta_2 \} \\
D_{\rm S} \vert E^{2k,\pm}_{\lambda,cl} &\cong& \bigoplus_{m_\lambda^\pm}
\bigl( \partial_x + {-{1\over 2}h' \pm \sqrt{\lambda} \over h} \bigr),\\
D_{\rm S} \vert E^{2k-1,\pm}_{\lambda,ccl} &\cong& \bigoplus_{m_\lambda^\pm}
\bigl( \partial_x + {{1\over 2}h' \pm \sqrt{\lambda} \over h} \bigr),
\end{eqnarray*}
where $m_\lambda^\pm := {\rm dim}\ E^{2k,\pm}_{\lambda,cl} =
{\rm dim}\ E^{2k-1,\pm}_{\lambda,ccl}$, the matrices $M_l$
are as in the table above,
$V := <\eta_1,\eta_2> \subset H_l$ for $1 \le l \le 4$, and
$\eta_1$ and $\eta_2$ are chosen as follows:

\begin{center}
\begin{tabular}{|c|c|c|} \hline
\rule[-1mm]{0mm}{7mm} $l$ & $\eta_1 \in$ & $\eta_2 =$ \\ \hline \hline
\rule[-5mm]{0mm}{12mm} $1$ & $E^{2k-2j-2}_{\lambda,cl}$
& ${1 \over \sqrt{\lambda}} *_N d_N \eta_1$ \\ \hline
\rule[-5mm]{0mm}{12mm} $2$ & $E^{2k-2j-3}_{\lambda,ccl}$ &
${1 \over \sqrt{\lambda}} d_N *_N \eta_1$ \\ \hline
\end{tabular} \qquad
\begin{tabular}{|c|c|c|} \hline
\rule[-1mm]{0mm}{7mm} $l$ & $\eta_1 \in$ & $\eta_2 =$ \\ \hline \hline
\rule[-5mm]{0mm}{12mm} $3$ & $E^{2k-2j-1}_{\lambda,cl}$
& ${1 \over \sqrt{\lambda}} *_N d_N \eta_1$ \\ \hline
\rule[-5mm]{0mm}{12mm} $4$ & $E^{2k-2j-2}_{\lambda,ccl}$ &
${1 \over \sqrt{\lambda}} d_N *_N \eta_1$ \\ \hline
\end{tabular} 
\end{center}
\end{lemma}

\bigskip
\centerline{\bf ACKNOWLEDGMENTS}

\medskip
Both authors are grateful to J. Br\"uning for many helpful discussions.
Moreover, the authors would like to express their gratitude to J. Cheeger,
R.T. Seeley, H. Moscovici, and the referee
for helpful comments on the subject.



\end{document}